\documentclass[12pt]{amsart}
\usepackage{graphicx}
\usepackage{amssymb,amstext,amsfonts,amscd}
\newtheorem{theorem}{Theorem}
\newtheorem{remark}[theorem]{Remark}
\newtheorem{lemma}[theorem]{Lemma}

\numberwithin{equation}{section}

\newcommand{\bn}{\mathbb{N}}
\newcommand{\bz}{\mathbb{Z}}
\newcommand{\br}{\mathbb{R}}
\newcommand{\bc}{\mathbb{C}}
\newcommand{\bp}{\mathbb{P}}

\newcommand{\bff}{\mathbb{F}}

\newcommand{\bee}{\mathbb{E}}
\newcommand{\bhh}{\mathbb{H}}

\newcommand{\rk}{\mathop\mathrm{rk}}

\newcommand{\MOD}{\ \mathrm{mod}\ }

\newcommand{\Proof}{\noindent{\it{Proof.}}\ \ }
\newcommand{\QED}{\ \ $\Box$}

\newcommand{\wvarphi}{{\widetilde{\varphi}}}
\newcommand{\wpsi}{{\widetilde{\psi}}}

\title{Real K3 surfaces with non-symplectic involution and applications. II}

\author{Viacheslav V. Nikulin and Sachiko Saito}

\address{Department of Pure Mathematics, The University of Liverpool,
Liverpool, L69 3BX, United Kingdom\\
Steklov Mathematical Institute, ul. Gubkina 8, Moscow, GSP-1,
Russia} \email{vnikulin@liv.ac.uk, vvnikulin@list.ru}

\address{
Department of Mathematics Education, Hakodate Campus, Hokkaido
University of Education, 1-2 Hachiman-cho, Hakodate 040-8567,
Japan} \email{sachi63@cc.hokkyodai.ac.jp}


\begin{document}
\pagestyle{plain}

\begin{abstract} We consider real forms of relatively minimal
rational surfaces $\bff_m$. Connected components of moduli of real
non-singular curves in $|-2K_{\bff_m}|$ had been classified
recently for $m=0,\,1,\,4$ in \cite{NikulinSaito05}. Applying
similar methods, here we fill the gap for $m=2$ and $m=3$ to
complete similar classification for any $0\le m\le 4$, when
$|-2K_{\bff_m}|$ is reduced.

The case of $\bff_2$ is especially remarkable and classical 
(quadratic cone in $\bp^3$). As an application, we finished
classification of connected components of moduli of real
hyper-elliptically polarized K3 surfaces and their deformations to
real polarized K3 surfaces started in \cite{NikulinSaito05} and
\cite{Nikulin05}. This could be important in some questions
because real hyper-elliptically polarized K3 surfaces can be
constructed explicitly.
\end{abstract}

\maketitle

\tableofcontents

\section{Introduction}\label{secIntro}

In \cite{NikulinSaito05} classification of connected components of
moduli of real K3 surfaces with non-symplectic involution had been
considered. As an application, classification of connected
components of moduli of real non-singular curves $A\in
|-2K_{\bff_m}|$ had been obtained for $m=0$, $1$, $4$ where
$\bff_m$ is a relatively minimal rational surface having
exceptional section with square $-m$. Here a K3 surface with
non-symplectic involution appears as the double coverings of
$\bff_m$ ramified in $A$ with the involution of the double
covering.

The main purpose of this paper is to fill the gap, and to get
similar classification for $m=2$ and $m=3$.

Case of $\bff_3$ (see Sect. \ref{secF3}) is very similar
to cases $m=0,\,1,\,4$
considered in \cite{NikulinSaito05} because in this case
we obtain K3 surfaces  with non-degenerate non-symplectic
involution.
This non-degenerate case had been treated in \cite{NikulinSaito05}
in general. Remarkably, enumeration of the connected components for
$\bff_3$ case is equivalent to classification of some special integral
polarized involutions classified in \cite{Nikulin79} in general.
Thus, direct application of results of \cite{Nikulin79} and
\cite{NikulinSaito05} solves the problem in this case.

The case of $\bff_3$ is related with some very classical problems.

Firstly, K3 surfaces of this case can be also considered as K3
surfaces with fixed elliptic pencil and its section such that the
elliptic pencil has no reducible fibres except one having the type
corresponding to the extended Dynkin diagram $\widetilde{\mathbb
A}_1$ or $\widetilde{\mathbb A}_2$. The inverse map of the group
law of the elliptic pencil with the section gives then the
non-symplectic involution. Thus, we classify connected components
of moduli of such real K3 surfaces. Similarly, the case of
$\bff_4$ which had been treated in \cite{NikulinSaito05} is
equivalent to K3 surfaces having fixed elliptic pencil with
section and no reducible fibres.

Secondly, the case of non-singular curves $A\in |-2K_{\bff_3}|$ is
equivalent to curves $A\in |-2K_{\bff_4}|$ having no singular
points except one quadratic (possibly degenerate) singular point.
Thus, we describe connected components of their moduli over $\br$
too. See Remark \ref{remF4F3}.

Thirdly, non-singular curves $A\in |-2K_{\bff_3}|$ are equivalent
to non-singular curves $A_1\in |-2K_Z|$ where $Z$ is a del Pezzo
surface having no singular points except one log-terminal singular
point of index two having the type $K_1$, moreover, the genus
$g(A_1)=9$ (see classification of log del Pezzo surfaces of index
$\le 2$ in \cite{AlexeevNikulin88}). Thus, we classify connected
components of moduli of real non-singular curves $A_1\in |-2K_Z|$
of such real del Pezzo surfaces $Z$. See Remark \ref{remZF3} for
details.

\medskip

The case of $\bff_2$ (see Sect. \ref{secF2}) when we consider
connected components of moduli of
real non-singular curves $A\in |-2K_{\bff_2}|$ is
especially remarkable. It gives K3 surfaces
with degenerate non-symplectic involution since $\bff_2$ has exceptional
curves with square $(-2)$. Thus, we
cannot directly apply results of \cite{NikulinSaito05} where only
non-degenerate case had been treated. Fortunately, all main
results of \cite{NikulinSaito05} survive in this case, and we
show that the connected component of moduli is also
determined by the action $(L,\tau,\varphi)$ of $\tau$ (holomorphic
non-symplectic involution) and $\varphi$ (anti-holomorphic
involution) on the homology lattice $L=H_2(X(\bc);\bz)$ of a K3 surface $X$
from this connected component. It is
clear that only in very special cases of degenerate
non-symplectic involutions on K3 similar result is valid. So,
we are very lucky. It would be interesting to describe all
similar nice cases in general, see Remark \ref{remrootinvcompF2}.
See \cite{AlexeevNikulin88} about classification of K3 surfaces over
$\bc$ with degenerate (i.e. arbitrary) non-symplectic involution.

Remarkably, enumeration of the connected components of moduli
in the case of $\bff_2$ is also equivalent to classification of
some special integral polarized involutions. Thus, again application
of results of \cite{Nikulin79} and \cite{NikulinSaito05} solves
the problem.

The case of $\bff_2$ is related with several classical problems.

Firstly, $\bff_2$ with contracted exceptional section is a quadratic cone in
$\bp^3$. Thus, we solve very classical problem of enumeration of
connected components of moduli of non-singular degree 8 curves on
the quadratic cone in 3-dimensional space (i. e. non-singular sections
of the cone by degree four forms). See Sect. \ref{subsecapp1F2}.

Deforming the cone to hyperboloid or ellipsoid, we obtain bidegree
(4,4) non-singular curves on hyperboloid or ellipsoid (the case of
$\bff_0$). Combining our results here and results of
\cite{NikulinSaito05} we completely classify these deformations.
This could be important in some very practical problems. See Sects
\ref{subsechypdefF2} and \ref{subsecelldefF2}.

Let us consider hyper-elliptically polarized K3 surfaces $(X,P)$
such that the linear system $|P|$ is hyper-elliptic and gives then
a double covering $|P|:X\to Y$. The case of $Y\cong \bff_2$
gives  {\it hyper-elliptically polarized K3 surfaces $(X,P)$ of
type $\bff_2$.} In this case, $|P|$ is double covering of
$\bff_2$ ramified in a non-singular curve $A\in |-2K_{\bff_2}|$.
Thus, applying our results, we obtain classification of connected
components of moduli of real hyper-elliptically polarized K3 surfaces
of type $\bff_2$. This finishes classification of connected components
of moduli of real hyper-elliptically polarized K3 surfaces since for
all other types similar classification had been obtained in
\cite{NikulinSaito05} (see also \cite{Nikulin05}). See Sect.
\ref{subsecconnhyp-ellpolK3F2}.

In Sect. \ref{subsecdefF2K3F2}, applying our results about
$\bff_2$ and results of \cite{NikulinSaito05} and
\cite{Nikulin05}, we classify deformations of real
hyper-elliptically polarized K3 surfaces of type $\bff_2$
to real polarized K3 surfaces. We completely describe all real
polarized K3 surfaces which can be obtained by such deformation.
Moreover, we enumerate all types of these deformations.
Similar classification for all other types of real
hyper-elliptically polarized K3 surfaces had been obtained in
\cite{NikulinSaito05} and \cite{Nikulin05}.
Thus, here we finalized these results.

All our results here and in \cite{NikulinSaito05},
\cite{Nikulin05} are based on Global Torelli Theorem \cite{PS71}
and surjectivity of Torelli map \cite{Kulikov77} for K3 surfaces,
and some geometric and arithmetic considerations. Thus, our
results are mainly existence and uniqueness results. It would be
very interesting to construct exactly all representatives of
connected components of moduli which we considered and classified.
Our results about deformations reduce everything to construction
of representatives of connected components of moduli of real
non-singular curves $A\in |-2K_{\bff_m}|$, $0\le m \le 4$, or
degree 6 curves in $\bp^2$. In many cases they are known by direct
construction. E. g. see \cite{Gudkov69}, \cite{Gudkov79},
\cite{GudkovShustin80}, \cite{Rokhlin78}, \cite{Viro80},
\cite{Zvonilov82}, \cite{Zvonilov83}, \cite{Zvonilov92}. It would
be very interesting to have such constructions for all connected
components of moduli. It will make our results effective and
useful for some questions which don't depend on deformation. This
construction is surely possible because all our connected
components of moduli are exactly labelled by some natural
invariants, and number of the connected components is bounded.

\section{Connected components of moduli of real non-singular
curves in $|-2K_{\bff_3}|$} \label{secF3}

We always consider algebraic varieties over $\bc$. If we say
``real", we consider algebraic varieties which are defined over
$\br$.

We denote by $\bff_m$, $m\ge 0$, a relatively minimal rational
surface with exceptional (if $m\ge 1$) section $s\cong \bp^1$
having $s^2=-m$. If $m=0$, then $\bff_0=\bp^1\times \bp^1$. If
$m>0$, then $\bff_m$ has a unique rational pencil with the section
$s$. We denote by $c$ a fibre of the rational pencil. The
canonical class $K_{\bff_m}=(-m-2)c-2s$ and
$-2K_{\bff_m}=(2m+4)c+4s$.

If $0\le m\le 2$, then a general element of the complete linear
system $|-2K_{\bff_m}|$ is a non-singular curve.

If $3\le m\le 4$, then $|-2K_{\bff_m}|=|(2m+4)c+3s|+s$ has the
fixed component $s$ where the general element $A_1\in
|(2m+4)c+3s|$ is non-singular. If $m=3$, then $A_1\cdot s=1$, and
$A=A_1+s\in|-2K_{\bff_3}|$ is union of two non-singular curves
which intersect transversally in one point (for $m=3$, we then
also call $A$ non-singular). If $m=4$, then $A_1\cdot s=0$, and
the whole curve $A=A_1+s\in |-2K_{\bff_4}|$ is non-singular.

If $m>4$, then $|-2K_{\bff_m}|=|(2m+4)c+2s|+2s$ has a non-reduced
fixed component $2s$.

Thus, cases $0\le m\le 4$ are characterized by the property that a
general element of the complete linear system $|-2K_{\bff_m}|$ is
reduced.

Connected components of moduli of real non-singular curves in
$|-2K_{\bff_m}|$ (for real $\bff_m$) were classified in
\cite{NikulinSaito05} for $m=0,\,1,\,4$. For $m=1$ it can be, in
principle, deduced from  \cite{Itenberg92}, \cite{Itenberg94}
where much more delicate classification had been considered. On
the other hand, our classification in \cite{NikulinSaito05} gives
{\it natural invariants} of connected components of moduli of real
non-singular curves in $|-2K_{\bff_1}|$.

In this paper, we want to complete this classification for reduced
$|-2K_{\bff_m}|$ (i. e. for $0\le m\le 4$) considering connected
components of moduli of real non-singular curves in
$|-2K_{\bff_3}|$ and $|-2K_{\bff_2}|$.

In this section, we consider the case of $|-2K_{\bff_3}|$ which is
very similar to $|-2K_{\bff_m}|$ for cases $m=0,\,1,\,4$
considered in  \cite{NikulinSaito05}. The case $|-2K_{\bff_3}|$ is
completely covered by the general theory which had been developed
in \cite{NikulinSaito05}. One needs only to perform necessary
calculations and geometric observations corresponding to this
case.

The case of $|-2K_{\bff_2}|$ is more delicate. It requires some
changes comparing to \cite{NikulinSaito05},  and it will be
considered in the next Sect. \ref{secF2}.

\subsection{Reduction to real K3 surfaces with
non-symplectic involution of type $(3,1,1)$}
\label{subsecK3F3}

Here we consider $\bff_3$. Like above,  $c$ denotes the fibre of
the rational pencil of $\bff_3$, and $s$ its exceptional section.
We have $c^2=0$, $s^2=-3$ and $c\cdot s=1$.

Let $A\in |-2K_{\bff_3}|$ be a non-singular curve, i. e. $A=A_1+s$
where $A_1$ is a non-singular curve of the genus $9$. Then
$A_1\cdot s=1$, and $A_1$ intersects $s$ transversally in one
point $A_1\cap s$. Let $f$ be the fibre of the rational pencil of
$\bff_3$ over this point.

Let $\sigma:Y\to \bff_3$ be the blow-up at the point $A_1\cap s$
and $e$ the exceptional curve of the blow-up. The $Y$ is a
rational surface with Picard number 3. We denote by $A_0$ the
proper pre-image of $s$, and by the same letters $A_1$, $f$ the
proper pre-images of $A_1$ and $f$ on $Y$. On $Y$, we have
$e^2=f^2=-1$, $A_0^2=-4$ and $A_0\cdot e=e\cdot f=1$, and
$A_0\cdot f=0$. Moreover, $A_1\cdot A_0=0$, $A_1\cdot e=1$ and
$A_1\cdot f=2$.

The proper pre-image $A=A_1+A_0\in |-2K_Y|$ of $A$ is a
non-singular curve. Let $\pi:X\to Y$ be the double covering of $Y$
ramified in $A\subset Y$. The $X$ is a K3 surface with the
non-symplectic involution $\tau$ of the double covering. Then
$A=X^\tau$ is identified with the fixed points set of the
involution, and $Y=X/\{1,\tau\}$ is the quotient space. We denote
$E=\pi^{-1}(e)$, $F=\pi^{-1}(f)$, and we identify $\pi^{-1}(A_0)$
with $A_0$. We have on $X$ that $A_0^2=E^2=F^2=-2$ and $A_0\cdot
E=1$, $A_0\cdot F=0$ and $E\cdot F=2$. The curves $A_0$ and $E$
are non-singular rational, i. e. they are isomorphic to $\bp^1$.
The curve $F$ is also non-singular rational if $A_1$ intersects
$f$ in two distinct points. If $A_1$ touches $f$, then
$F=F^\prime+F^{\prime\prime}$ is union of two non-singular
rational curves which are conjugate by $\tau$, and $F^\prime\cdot
F^{\prime\prime}=1$.

The determinant of the Gram matrix of $A_0$, $E$ and $F$ is $2$.
It follows that $A_0$, $E$ and $F$ generate a primitive
3-dimensional sublattice $S$ in the Picard lattice $S_X$ of $X$
which is the eigenvalue $1$ sublattice
$S=S_X^\tau=H_2(X(\bc);\bz)^\tau$ of the involution $\tau$ (since
$Y$ has the Picard number 3). Thus, the lattice $S=S_X^\tau$ is a
hyperbolic lattice of the rank $r(S)=3$. The discriminant group
$S^\ast/S\cong (\bz/2\bz)^1$ is a 2-elementary group of the rank
$a(S)=1$. The element $cl(F)/2\in S^\ast$ has $(F/2)^2=-2/4\not\in
\bz$. Thus, the 2-elementary lattice $S$ has the invariant
$\delta(S)=1$ (it would be $0$ if any element $x^\ast\in S^\ast$
has $(x^\ast)^2\in \bz$). Thus, finally, the 2-elementary
hyperbolic lattice $S$ has the invariants
$(r(S),a(S),\delta(S))=(3,1,1)$ which are the main invariants of
the K3 surface with non-symplectic involution $(X,\tau)$.

Vice versa, by \cite{AlexeevNikulin88}, any K3 surface with
non-symplectic involution $(X,\tau)$ having the lattice
$S=S_X^\tau$ with the same invariants
$(r(S),a(S),\delta(S))=(3,1,1)$ can be obtained in the same way as
above from a non-singular curve $A_1+s\in |-2K_{\bff_3}|$.
Moreover, all exceptional curves (i. e. irreducible with negative
self-intersection) on $Y=X/\{ 1,\tau\}$ are exactly the curves
$A_0$, $e$ and $f$. Thus, one obtains the non-singular pair
$A_1+s\in |-2K_{\bff_3}|$ from $(X,\tau)$ in a unique way by
contracting $e$ and taking the image of $A=X^\tau$ which consists
of disjoint union of a curve $A_1$ of genus 9 and $A_0\cong
\bp^1$. Thus, pairs $(X,\tau)$ with the invariants $(3,1,1)$ and
non-singular $A_1+s\in |-2K_{\bff_3}|$ are in one to one
correspondence. In particular, they have the same moduli space and
its connected components.

We mention that K3 surfaces with non-symplectic involution
$(X,\tau)$ of type $(3,1,1)$ give a very important class of K3
surfaces. They are K3 surfaces with elliptic pencil $|E+F|$ (it is
pre-image of the rational pencil of $\bff_3$), its fixed section
$A_0$ and a reducible fibre $E+F$ satisfying conditions (a), (b)
and (c) from \eqref{ellfibrtypeF3} below:
\begin{equation}
\label{ellfibrtypeF3}
\end{equation}
\noindent (a) $E$ is irreducible non-singular rational and $E\cdot
A_0=1$;

\noindent (b) $F$ is either irreducible non-singular rational or
is sum of two irreducible non-singular rational curves, and
$E\cdot F=2$ (i. e. the degenerate fibre corresponds to the
extended root system $\widetilde{\mathbb{A}}_1$ or
$\widetilde{\mathbb{A}}_2$);

\noindent (c) the elliptic pencil does not have other reducible
fibres.

\medskip

\noindent Then $\tau$ is the inverse map of the group law of the
elliptic pencil with the section. It is a non-symplectic
involution of type $(3,1,1)$.

\medskip

Similarly, the case of non-singular curves $A\in |-2K_{\bff_4}|$
which had been treated in \cite{NikulinSaito05} is equivalent to K3
surfaces with fixed elliptic pencil with section and without reducible
fibres (in this case, $(r(S),a(S),\delta(S))=(2,0,0)$).

\medskip

Now let us assume that $\bff_3$ and $A=A_1+s\in |-2K_{\bff_3}|$
are real. Then the point $A_1 \cap s$ is real and $\bff_3(\br)$ is
non-empty. Then $Y(\br)$ is non-empty and the corresponding
anti-holomorphic involution $\theta$ of $Y$ can be lifted to an
anti-holomorphic involution $\varphi$ of the K3 surface $X$. Thus,
we obtain a real K3 surface with non-symplectic involution
$(X,\tau, \varphi)$, i. e. $\tau$ and $\varphi$ commute. The lift
$\varphi$ is unique up to replacement by the {\it related
anti-holomorphic involution} $\widetilde{\varphi}=\varphi\tau$.
This choice is defined by the image
$\pi(X_\varphi(\br))=A^+(\varphi)\subset Y(\br)$ where
$X_\varphi(\br)=X(\bc)^\varphi$ is the set of real points of the
real K3 surface $(X,\varphi)$. Here $A^+(\varphi)$ is called the
{\it positive curve of $A$.} If one takes the related real
structure $\widetilde{\varphi}=\varphi\tau$ on $X$, then the
positive curve becomes
$\pi(X_{\widetilde{\varphi}}(\br))=A^+(\widetilde{\varphi})=
A^-(\varphi)=\overline{Y(\br)-A^+(\varphi)}$. We have
$A^+(\varphi)\cap A^+(\widetilde{\varphi})=A(\br)$.

Thus, classifying connected components of moduli of real K3
surfaces with non-symplectic involution $(X,\tau,\varphi)$ of type
$(3,1,1)$, we actually solve more delicate problem of
classification of connected components of moduli of positive real
non-singular curves $A^+$ where $A\in |-2K_{\bff_3}|$ is
non-singular. Here $A^+$ is defined by a choice of one of two
parts of $\bff_3(\br)-A(\br)$ where $A^+\cap
(\bff_3(\br)-A^+)=A(\br)$.

\subsection{Enumeration of connected components of moduli of
real K3 surfaces with non-symplectic involution of type
$(3,1,1)$} \label{subsecModinv} Here we describe connected
components of moduli of real K3 surfaces with non-symplectic
involution of type $(3,1,1)$.

For any $(X,\tau)$ of type $(3,1,1)$,  the surface $Y$ has no
exceptional curves with square $(-2)$. Thus, $(X,\tau)$ is always
non-degenerate in the sense of (\cite{NikulinSaito05}, Sect. 2.1).
In this case, a connected component of moduli of real K3 surfaces
with non-symplectic involution $(X,\tau,\varphi)$ is determined by
the isomorphism class of the action $(L,\tau, \varphi)$ of $\tau$
and $\varphi$ on the homology lattice $L=H_2(X(\bc);\bz)$ with
intersection pairing (\cite{NikulinSaito05}, Theorem 1). Here $L$
is an even unimodular lattice of signature $(3,19)$.

Since $A_0$, $e$ and $f$ are all exceptional curves on $Y$ and
they  have different geometry on $Y$, the anti-holomorphic
involution $\theta$ on $Y$ should send each of these three curves
to itself changing its orientation. It follows that
$\varphi\,|\,S=\varphi\,|\,L^\tau=-1$. Moreover, the eigenvalue
$1$ part $L^\varphi$ of $\varphi$ in $L$ must be hyperbolic (it
has exactly one positive square). By (\cite{NikulinSaito05},
Theorem 1), any isomorphism class of triplets
\begin{equation}
(L,\tau,\varphi) \label{triplinvF3}
\end{equation}
with these properties is possible. These isomorphism classes are
in one-to-one correspondence with connected components of moduli.
Moreover by (\cite{NikulinSaito05}, Theorem 14 and Proposition 15)
the isomorphism classes are determined by their genus invariants.

The genus invariants can be enumerated using general results of
\cite{Nikulin83} which are valid for any $S$ (see also their
adaptation to the particular case of K3 with non-symplectic
involution in \cite{NikulinSaito05}). On the other hand, our case
is very simple, and it can be also reduced to invariants of
integral polarized involutions (see Sect. 3 in \cite{Nikulin79}).
In this case, all necessary calculations are already done in
Theorem 3.4.3 of \cite{Nikulin79}. We prefer the second
opportunity.

The lattice $S=U \oplus \bz cl(F)$ where $U=\bz cl(A_0)+\bz
(cl(E)+cl(F))$ is the hyperbolic even unimodular lattice of
signature $(1,1)$, and $\bz cl(F)$ is isomorphic to $\langle -2
\rangle$ since $F^2=-2$. Here $\oplus$ means orthogonal sum of
lattices, and $\langle M\rangle$ a lattice with the symmetric
matrix $M$. Since $\tau\,|\,U=1$ and $\varphi\,|\,U=-1$, a triplet
$(L,\tau,\varphi)$ defines an integral polarized involution
\begin{equation}
\label{intpolinvF3}
 (L_1,\varphi_1,cl(F)).
\end{equation}
Here $L_1$ is the orthogonal complement to $U$ in $L$ (thus, it is
an even unimodular lattice of signature $(2,18)$,
$\varphi_1=\varphi\,|\,L_1$ and $L_1^{\varphi_1}$ is hyperbolic,
moreover, $\varphi_1(cl(F))=-cl(F)$ and $F^2=-2$.

Taking orthogonal sum of \eqref{intpolinvF3} with $U$, one obtains
the triplet \eqref{triplinvF3}. Thus, their isomorphism classes
are in one to one correspondence.


Genus invariants of integral polarized involutions
\eqref{intpolinvF3} (and then \eqref{triplinvF3}) were completely
classified in (\cite{Nikulin79}, Theorem 3.4.3). Thus, they
enumerate the connected components of moduli. Below we present
these results. (One should multiply by $(-1)$ the pairing in $L_1$
and apply Theorem 3.4.3 in \cite{Nikulin79} to $l_{(+)}=18$,
$l_{(-)}=2$, $t_{(-)}=1$, $t_{(+)}=r-1$ and $n=2$.)

The complete genus and isomorphism class invariants of
\eqref{intpolinvF3} (and then \eqref{triplinvF3}) are
\begin{equation}
(r,a,\delta_\varphi;\delta_F,\delta_{\varphi F}).
 \label{geninvF3}
\end{equation}
Here $r=\rk L^\varphi\in \bn$;
$((L^\varphi)^\ast/L^\varphi)=(\bz/2\bz)^a$ where $a\ge 0$ is an
integer; $\delta_\varphi \in \{0,1\}$ is equal to $0$ if and only
if $x\cdot \varphi(x)\equiv 0\mod 2$ for any $x\in L$. They are
complete invariants of the corresponding pair $(L,\varphi)$. Here
$\delta_{F}\in \{0,1\}$ is equal to $0$ if and only if $F \cdot
L_\varphi\equiv 0\mod 2$; here $\delta_{\varphi F}\in \{0,1\}$ is
equal to $0$ if and only if $x\cdot \varphi(x)\equiv x\cdot F$ for
any $x\in L$. Here and in what follows, we always denote by
$L^\varphi$ and $L_\varphi$ the eigenvalue $1$ and $-1$ parts
respectively of the action of an involution $\varphi$ on a module
$L$.

The invariants \eqref{geninvF3} must satisfy the following
relations which are sufficient and necessary for existence of a
triplet \eqref{intpolinvF3} or the corresponding triplet
\eqref{triplinvF3} with these invariants. Thus, these relations
enumerate connected components of moduli.

\noindent {\bf 0. Conditions on $(r,a,\delta_\varphi )$:}

\noindent (1) $1\le r\le 18$, $0\le a\le \min \{r,\ 20-r\}$;

\noindent (2) $r+a\equiv 0\mod 2$; if $\delta_\varphi=0$, then
$r\equiv 2\mod 4$;

\noindent (3) if $a=0$, then $(\delta_\varphi=0,\ r\equiv 2\mod
8)$;

\noindent (4) if $a=1$, then $r\equiv 1,\,3\mod 8$;

\noindent (5) if $(a=2,\  r\equiv 6\mod 8)$, then
$\delta_\varphi=0$;

\noindent (6) if $(a=r,\ \delta_\varphi=0)$, then $r\equiv 2\mod
8$;

\noindent (7) if $(a=20-r, \delta_\varphi=0)$, then $r\equiv 2\mod
8$.

\medskip

\noindent{\bf I. Conditions on $\delta_F$, $\delta_{\varphi
 F}$:}

\noindent {\bf General conditions:}

\noindent (1) if $\delta_F=0$, then $\delta_{\varphi}=1$;

\noindent (2) if $\delta_{\varphi F}=0$, then $(\delta_F=0,\
\delta_{\varphi }=1,\ r\equiv 1 \mod 4)$.

\medskip

\noindent {\bf Relations near the boundary $a=20-r$:}

\noindent (3) if $a=20-r$, then $\delta_F=0$;

\noindent (4) if $(a=20-r,\ \delta_{\varphi F}=0)$, then $r\equiv
1\mod 8$.

\medskip

\noindent {\bf Relations near the boundary $a=0$:}

\noindent (5) if $a=0$, then $\delta_F=1$;

\noindent (6) if $(a=1,\ \delta_F=0)$, then $(\delta_{\varphi
F}=0,\ r\equiv 1\mod 8)$;

\noindent (7) if $(a=2,\ \delta_F=0)$, then $r\equiv 0,\, 2\mod
8$;

\noindent (8) if $(a=3,\ \delta_F=0,\ r\equiv 5 \mod 8)$, then
$\delta_{\varphi F}=0$.

\medskip

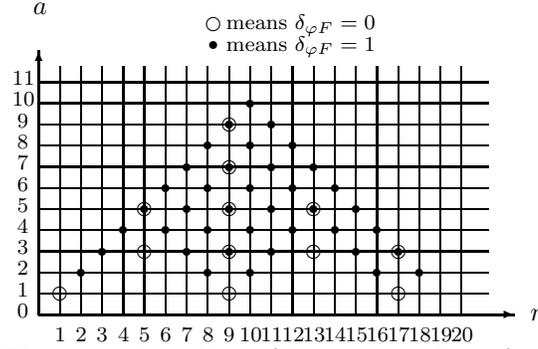
\begin{figure}
\begin{picture}(200,140)
\put(66,110){\circle{5}} \put(71,108){{\tiny means
$\delta_{\varphi F}=0$}} \put(66,102){\circle*{3}}
\put(71,100){{\tiny means $\delta_{\varphi F}=1$}}

\multiput(8,0)(8,0){20}{\line(0,1){94}}
\multiput(0,8)(0,8){11}{\line(1,0){170}}
\put(0,0){\vector(0,1){100}} \put(0,0){\vector(1,0){180}} \put(
6,-10){{\tiny $1$}} \put( 14,-10){{\tiny $2$}} \put(
22,-10){{\tiny $3$}} \put( 30,-10){{\tiny $4$}} \put(
38,-10){{\tiny $5$}} \put( 46,-10){{\tiny $6$}} \put(
54,-10){{\tiny $7$}} \put( 62,-10){{\tiny $8$}} \put(
70,-10){{\tiny $9$}} \put( 76,-10){{\tiny $10$}} \put(
84,-10){{\tiny $11$}} \put( 92,-10){{\tiny $12$}}
\put(100,-10){{\tiny $13$}} \put(108,-10){{\tiny $14$}}
\put(116,-10){{\tiny $15$}} \put(124,-10){{\tiny $16$}}
\put(132,-10){{\tiny $17$}} \put(140,-10){{\tiny $18$}}
\put(148,-10){{\tiny $19$}} \put(156,-10){{\tiny $20$}}

\put(-8, -1){{\tiny $0$}} \put(-8,  7){{\tiny $1$}} \put(-8,
15){{\tiny $2$}} \put(-8, 23){{\tiny $3$}} \put(-8, 31){{\tiny
$4$}} \put(-8, 39){{\tiny $5$}} \put(-8, 47){{\tiny $6$}} \put(-8,
55){{\tiny $7$}} \put(-8, 63){{\tiny $8$}} \put(-8, 71){{\tiny
$9$}} \put(-10, 79){{\tiny $10$}} \put(-10, 87){{\tiny $11$}}
\put( -2,114){{\footnotesize $a$}} 
\put(186, -2){{\footnotesize $r$}} 


\put(  8, 8){\circle{5}}
\put( 72, 8){\circle{5}} \put( 72, 24){\circle{5}} \put( 72,
40){\circle{5}} \put( 72, 56){\circle{5}} \put( 72,
72){\circle{5}}
\put(136, 8){\circle{5}} \put(136, 24){\circle{5}}
\put( 16,16){\circle*{3}}
\put( 64,16){\circle*{3}} \put( 80,16){\circle*{3}}
\put(128,16){\circle*{3}} \put(144,16){\circle*{3}} \put(
24,24){\circle*{3}}
\put( 40,24){\circle{5}} \put( 40,40){\circle{5}} \put(
56,24){\circle*{3}} \put( 72,24){\circle*{3}} \put(
88,24){\circle*{3}}
\put(104,24){\circle{5}} \put(104,40){\circle{5}}
\put(120,24){\circle*{3}}
\put(136,24){\circle*{3}} \put( 32,32){\circle*{3}} \put(
48,32){\circle*{3}} \put( 64,32){\circle*{3}} \put(
80,32){\circle*{3}} \put( 96,32){\circle*{3}}
\put(112,32){\circle*{3}} \put(128,32){\circle*{3}} \put(
40,40){\circle*{3}} \put( 56,40){\circle*{3}} \put(
72,40){\circle*{3}} \put( 88,40){\circle*{3}}
\put(104,40){\circle*{3}} \put(120,40){\circle*{3}} \put(
48,48){\circle*{3}} \put( 64,48){\circle*{3}} \put(
80,48){\circle*{3}} \put( 96,48){\circle*{3}}
\put(112,48){\circle*{3}} \put( 56,56){\circle*{3}} \put(
72,56){\circle*{3}} \put( 88,56){\circle*{3}}
\put(104,56){\circle*{3}} \put( 64,64){\circle*{3}} \put(
80,64){\circle*{3}} \put( 96,64){\circle*{3}} \put(
72,72){\circle*{3}} \put( 88,72){\circle*{3}}
\put( 80,80){\circle*{3}}   
\end{picture}

\caption{$\bff_3$: The complete list of invariants \ref{geninvF3}
of connected components of moduli with $\delta_F=0$ (then
$\delta_{\varphi}=1$).} \label{diag0F3}
\end{figure}

\begin{figure}
\begin{picture}(200,140)
\put(66,110){\circle{7}} \put(71,108){{\tiny means
$\delta_\varphi=0$,}} \put(66,102){\circle*{3}}
\put(71,100){{\tiny means $\delta_\varphi=1$}}

\multiput(8,0)(8,0){20}{\line(0,1){94}}
\multiput(0,8)(0,8){11}{\line(1,0){170}}
\put(0,0){\vector(0,1){100}} \put(0,0){\vector(1,0){180}} \put(
6,-10){{\tiny $1$}} \put( 14,-10){{\tiny $2$}} \put(
22,-10){{\tiny $3$}} \put( 30,-10){{\tiny $4$}} \put(
38,-10){{\tiny $5$}} \put( 46,-10){{\tiny $6$}} \put(
54,-10){{\tiny $7$}} \put( 62,-10){{\tiny $8$}} \put(
70,-10){{\tiny $9$}} \put( 76,-10){{\tiny $10$}} \put(
84,-10){{\tiny $11$}} \put( 92,-10){{\tiny $12$}}
\put(100,-10){{\tiny $13$}} \put(108,-10){{\tiny $14$}}
\put(116,-10){{\tiny $15$}} \put(124,-10){{\tiny $16$}}
\put(132,-10){{\tiny $17$}} \put(140,-10){{\tiny $18$}}
\put(148,-10){{\tiny $19$}} \put(156,-10){{\tiny $20$}}

\put(-8, -1){{\tiny $0$}} \put(-8,  7){{\tiny $1$}} \put(-8,
15){{\tiny $2$}} \put(-8, 23){{\tiny $3$}} \put(-8, 31){{\tiny
$4$}} \put(-8, 39){{\tiny $5$}} \put(-8, 47){{\tiny $6$}} \put(-8,
55){{\tiny $7$}} \put(-8, 63){{\tiny $8$}} \put(-8, 71){{\tiny
$9$}} \put(-10, 79){{\tiny $10$}} \put(-10, 87){{\tiny $11$}}
\put( -2,114){{\footnotesize $a$}} 
\put(186, -2){{\footnotesize $r$}} 

\put( 16, 0){\circle{7}} \put( 80, 0){\circle{7}} \put(144,
0){\circle{7}} \put( 16,16){\circle{7}} \put( 48,16){\circle{7}}
\put( 80,16){\circle{7}} \put(112,16){\circle{7}}
\put( 48,32){\circle{7}} \put(
80,32){\circle{7}} \put(112,32){\circle{7}} \put(
80,48){\circle{7}}
\put(
80,64){\circle{7}}

\put(  8, 8){\circle*{3}} \put( 24, 8){\circle*{3}} \put( 72,
8){\circle*{3}} \put( 88, 8){\circle*{3}} \put(136,
8){\circle*{3}}
\put(
16,16){\circle*{3}} \put( 32,16){\circle*{3}} \put(
64,16){\circle*{3}} \put( 80,16){\circle*{3}} \put(
96,16){\circle*{3}} \put(128,16){\circle*{3}}
\put( 24,24){\circle*{3}} \put(
40,24){\circle*{3}} \put( 56,24){\circle*{3}} \put(
72,24){\circle*{3}} \put( 88,24){\circle*{3}}
\put(104,24){\circle*{3}} \put(120,24){\circle*{3}}
\put( 32,32){\circle*{3}} \put(
48,32){\circle*{3}} \put( 64,32){\circle*{3}} \put(
80,32){\circle*{3}} \put( 96,32){\circle*{3}}
\put(112,32){\circle*{3}}
\put(40,40){\circle*{3}} \put( 56,40){\circle*{3}} \put(
72,40){\circle*{3}} \put( 88,40){\circle*{3}}
\put(104,40){\circle*{3}}
\put(48,48){\circle*{3}} \put( 64,48){\circle*{3}} \put(
80,48){\circle*{3}} \put( 96,48){\circle*{3}}
\put( 56,56){\circle*{3}}
\put(72,56){\circle*{3}} \put( 88,56){\circle*{3}}
\put( 64,64){\circle*{3}}
\put(80,64){\circle*{3}}
\put(
72,72){\circle*{3}}
\end{picture}

\caption{$\bff_3$: The complete list of invariants \ref{geninvF3}
of connected components of moduli with $\delta_F=1$ (then
$\delta_{\varphi F}=1$).} \label{diag1F3}
\end{figure}
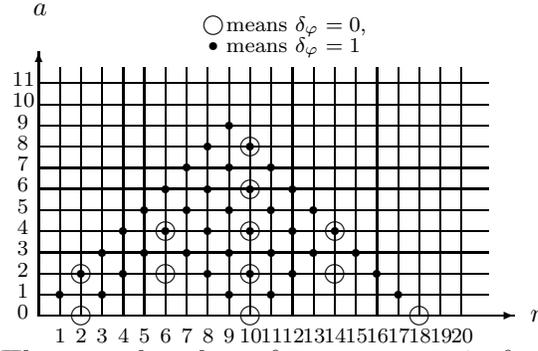

It is easy to enumerate all invariants \eqref{geninvF3} satisfying
these conditions. They are all presented in Figure \ref{diag0F3}
for $\delta_F=0$, and in Figure \ref{diag1F3} for $\delta_F=1$.
Thus Figures \ref{diag0F3} and \ref{diag1F3} enumerate all
connected components of moduli.


\medskip

In \cite{NikulinSaito05}, instead of $\delta_F$,
$\delta_{\varphi}$ and $\delta_{\varphi F}$, other invariants were
also used. They are: the subgroup $H\subset
(\bz/2\bz)cl(F)=(2S^\ast/2S)$, the invariant $\delta_{\varphi
S}\in \{0,1\}$, the characteristic element $v\in H$ when
$\delta_{\varphi S}=0$ (the $v$ is not defined when
$\delta_{\varphi S}=1$). They are related in the obvious way:
$H=0$ iff $\delta_F=1$; $\delta_{\varphi S}=0$ and $v=0$ ($\mod
2S$) iff $\delta_{\varphi}=0$; $\delta_{\varphi S}=0$ and $v\sim
F$ $(\mod 2S)$ iff $\delta_{\varphi F}=0$; otherwise (if
$\delta_{\varphi}=\delta_{\varphi F}=1$) the invariant
$\delta_{\varphi S}=1$.

By simple direct calculations or applying (\cite{NikulinSaito05},
Theorem 11 and Corollary 13), one finds how invariants
\eqref{geninvF3} of related involutions $\varphi$ and
$\widetilde{\varphi}=\tau\varphi$  can be calculated from each
other:
\begin{equation}
\begin{array}{l}
r(\varphi)+r(\widetilde{\varphi})=19;\
a({\varphi})+\delta_F(\varphi)=
a(\widetilde{\varphi})+\delta_F(\widetilde{\varphi});\\
\delta_F(\varphi)+\delta_F(\widetilde{\varphi})=1; \
\delta_{\widetilde{\varphi}}=\delta_{\varphi F};\
\delta_{\widetilde{\varphi} F}=\delta_\varphi\,.
\end{array}
\label{relatedF3}
\end{equation}

One can clearly see all these relations looking at Figures
\ref{diag0F3}, \ref{diag1F3}.

Thus, we finally obtain the main result which enumerates connected
components of moduli we are looking for.

\begin{theorem} There are exactly 102 connected components of
moduli of positive real non-singular curves $A\in |-2K_{\bff_3}|$
up to the action of $Aut(\bff_3/\br)$, equivalently of real K3
surfaces with non-symplectic involution of type $(3,1,1)$,
equivalently of real elliptic K3 surfaces with section and no
reducible fibres except one of the type $\widetilde{\mathbb{A}}_1$
or $\widetilde{\mathbb{A}}_2$. They are labelled by their natural
invariants \eqref{geninvF3} satisfying relations 0.1--7 and I.1--8
above enumerated in Figures \ref{diag0F3} and \ref{diag1F3}.

Identifying related pairs of positive curves and anti-holomorphic
involutions of K3, there are exactly 51 connected components of
moduli of real non-singular curves $A\in |-2K_{\bff_3}|$ up to the
action of $Aut(\bff_3/\br)$, equivalently of related pairs
$(X,\tau, \varphi)$ and $(X,\tau,\widetilde{\varphi}=\tau\varphi)$
of real K3 surfaces with non-symplectic involution of type
$(3,1,1)$, equivalently of related by the inverse map pairs of
real elliptic K3 surfaces with section and no reducible fibres
except one of the type $\widetilde{\mathbb{A}}_1$ or
$\widetilde{\mathbb{A}}_2$. They are labelled by their related by
\eqref{relatedF3} pairs of natural invariants \eqref{geninvF3}
satisfying relations 0.1--7 and I.1--8 above enumerated in Figures
\ref{diag0F3} and \ref{diag1F3}.

\label{themainF3}
\end{theorem}

Below we consider geometric interpretation of the invariants
\eqref{geninvF3} in terms of K3 and curves on $\bff_3$.

\subsection{Geometric interpretation of the invariants \eqref{geninvF3}
in terms of K3} \label{geomK3F3} We have the following geometric
interpretation of the invariants
$(r,a,\delta_\varphi,\delta_{\varphi F})$ in terms of the real K3
surface $(X,\varphi)$ (see \cite{Nikulin79} and also
\cite{Kharlamov75a}). We denote by $T_g$ an orientable closed
surface of genus $g$. We have
\begin{equation}
X_\varphi(\br)
= \left\{
\begin{array}{ll}
\emptyset         & \text{if $(r,a,\delta_\varphi)=(10,10,0)$}\\
T_1 \amalg T_1    & \text{if $(r,a,\delta_\varphi)=(10,8,0)$}\\
T_g \amalg k T_0  & \text{otherwise,}
\end{array}
\right. \label{realcomponents}
\end{equation}
where
$$
g=\frac{22-r-a}{2},\ \ \ k=\frac{r-a}{2};
$$
\begin{equation}
X_\varphi(\br) \sim 0 \mod 2\ \ \ \text{\ in\ } H_2(X(\bc);\bz),
\label{geomdeltaphi}
\end{equation}
if and only if $\delta_\varphi=0$, and
\begin{equation}
X_\varphi(\br) \sim F(\bc)  \mod 2\ \ \ \text{\ in\ }
H_2(X(\bc);\bz), \label{geomdeltaphiF}
\end{equation}
if and only if $\delta_{\varphi F}=0$.

For the invariant $\delta_{F}$, we have the following theorem (see
Theorem 16 in \cite{NikulinSaito05}).

\begin{theorem}[\cite{NikulinSaito05}]
\label{realcurve} Let $X$ be a compact K\"{a}lerian surface with
an anti-holomorphic involution $\varphi$, and
$X_\varphi(\br)=X(\bc)^\varphi$ the set of real points. We assume
that $H_1(X(\bc);\bz) = 0$ and $X_\varphi(\br) \neq \emptyset$.
Let $C$ be a $1$-dimensional complex submanifold of $X$ with
$\varphi(C)=C$, and $C(\br)$ be the fixed point set of $\varphi$
on $C$. Let $cl(C)\in H_2(X(\bc);\bz)$ and $cl(C(\br))\in
H_1(X_\varphi(\br);\bz/2\bz)$ denote the homology classes
represented by $C$ and $C(\br)$ respectively.

Then, $cl(C)\cdot H_2(X(\bc);\bz)_\varphi\equiv 0\mod 2$
(equivalently, $\delta_C(\varphi)=0$) if and only if $cl(C(\br)) =
0$ in $H_1(X_\varphi(\br);\bz/2\bz)$ (as usual,
$H_2(X(\bc);\bz)_\varphi$ denotes the eigenvalue $-1$ part of
$\varphi$).
\end{theorem}

\bigskip

Thus, if $F$ is irreducible, then $\delta_F=0$ if and only if
$F(\br )=0$ in $H_1(X_\varphi(\br);\bz/2\bz)$. Changing the real
structure a little (in the same connected component of moduli), we
can always make $F$ irreducible, $F\cong \bp^1$, and it is then a
complex submanifold of $X(\bc)$. By deformation argument, we then
have the same statement for a reducible $F$. Since
$X_\varphi(\br)$ is always non-empty, we get at any case
\begin{equation}
\delta_F=0\ if\ and\ only\ if\ cl(F(\br))=0\ in\
H_1(X_\varphi(\br);\bz/2\bz). \label{deltaFF3}
\end{equation}

Thus, all invariants \eqref{geninvF3} have very clear geometric
meaning in terms of the real K3 surface $(X,\varphi)$.

\begin{theorem} The connected component of moduli of a real K3 surface with
non-symplectic involution $(X,\tau,\varphi)$ of the type $(3,1,1)$
(equivalently of a positive real non-singular curve $A\in
|-2K_{\bff_3}|$) is determined by the geometric invariants:
topology of the real part $X_\varphi(\br)$, if $F_\varphi(\br)\sim
0\mod 2$ in $H_1(X_\varphi(\br); \bz)$, if $X_\varphi (\br)\sim
0\mod 2$ in $H_2(X(\bc);\bz)$, if $X_\varphi (\br)\sim F(\bc)\mod
2$ in $H_2(X(\bc);\bz)$.

See \eqref{realcomponents}, \eqref{geomdeltaphi},
\eqref{geomdeltaphiF} and \eqref{deltaFF3} about their exact
relation with the complete invariants \eqref{geninvF3} of the
connected component of moduli. \label{thegeomK31F3}
\end{theorem}

Taking the related involution $\widetilde{\varphi}=\tau\varphi$
instead of $\varphi$, we obtain geometric interpretation of the
invariants \eqref{geninvF3} in terms of the real K3 surface $(X,
{\widetilde{\varphi}})$.

\begin{theorem} The connected component of moduli of a real K3 surface with
non-symplectic involution $(X,\tau,\varphi)$ of type $(3,1,1)$
(equivalently of a positive real non-singular curve $A\in
|-2K_{\bff_3}|$) is determined by the geometric invariants:
topology of the real part $X_{\widetilde{\varphi}}(\br)$, if
$F_{\widetilde{\varphi}}(\br)\sim 0\mod 2$ in
$H_1(X_{\widetilde{\varphi}}(\br); \bz)$, if
$X_{\widetilde{\varphi}} (\br)\sim 0\mod 2$ in $H_2(X(\bc);\bz)$,
if $X_{\widetilde{\varphi}} (\br)\sim F(\bc)\mod 2$ in
$H_2(X(\bc);\bz)$.

See \eqref{realcomponents}, \eqref{geomdeltaphi},
\eqref{geomdeltaphiF} \eqref{deltaFF3} and \eqref{relatedF3} about
their exact relation with the complete invariants \eqref{geninvF3}
of the connected component of moduli. \label{thegeomK32F3}
\end{theorem}

\subsection{Geometric interpretation of the invariants \eqref{geninvF3}
in terms of $A\in |-2K_{\bff_3}|$} \label{geomF3F3} Next we
consider the dividedness of the fixed point curve of a real K3
surface with a non-symplectic involution.

Let $A$ be a possibly disconnected non-singular complex compact
curve with an anti-holomorphic involution $\varphi$, and $A(\br)$
the fixed point set of $\varphi$. We remind that $(A,\varphi)$ (or
simply, $A$) is called {\it dividing} if $A(\br)=0\ \in
H_1(A;\bz/2\bz)$. We have the following result which had been
partially used in \cite{NikulinSaito05} (see the proof of Theorem
26 in \cite{NikulinSaito05}).

\begin{theorem}\label{div}
Let $(X,\tau,\varphi)$ be a real K3 surface with non-symplectic
involution of type $(S,\theta)=(H_2(X(\bc);\bz)^\tau,\, \varphi |
H_2(X(\bc);\bz)^\tau )$ (see \cite{NikulinSaito05}); let
$A=X^\tau$ be the fixed point set of $\tau$ in $X$. Assume that
$A$ is non-empty. Then, the real curve $A$ is dividing if and only
if $\delta_{\varphi S}=0$ (equivalently, $A(\br)\mod 2\in
H_2(X(\bc);\bz)^\tau\mod 2$).
\end{theorem}

\begin{proof}
The main method is a so-called Donaldson's trick (see \cite{D} or
\cite{DIK2000}). It is known that the K3 surface $X$ is algebraic.
Hence, we can take a hyperplane section class $h$ in $S$. We
denote $L=H_2(X(\bc);\bz)$. Since $\varphi_*(h)=-h$, $h$ is
contained in $S_- = L^{\tau} \cap L_{\varphi}$. Moreover, $h$ is
contained in the K\"ahler cone $\mathcal{C}^+_X$ of $X$. (See
\cite{BPV}, p.242, for example.) Thus we see $S_- \cap
\mathcal{C}^+_X \neq \emptyset$. Let $I$ be the complex structure
of the K3 surface $X$ and take a K\"ahler class $c$ of $(X,I)$ in
$((L^{\tau} \cap L_{\varphi}) \otimes \br) \cap \mathcal{C}^+_X$.
Since $c_1(X)=0$, by the Calabi-Yau theorem (for example,
\cite{Besse}, Theorem 11.15) $0$ is the Ricci form of one and only
one K\"ahler metric in the class $c$. Thus we can take the unique
K\"ahler form $P$ with zero Ricci form in the class $c$, where
$c=[P]$. Let $g$ be the K\"ahler metric whose K\"ahler form is
$P$, i.e., $P=g(I(\ \cdot\ ),\ \cdot  )$. By the uniqueness of the
Ricci flat K\"ahler form, we have $\tau^*P = P$ and $\varphi^*P =
- P$.
Since $\tau$ is holomorphic on $(X,I)$, we have $\tau^*g=g$. And,
since $\varphi$ is anti-holomorphic on $(X,I)$, we have $\varphi^*g=g$.

Now we can take a nowhere vanishing holomorphic $2$-form
$\omega_I$ on $(X,I)$
such that $\tau^* \omega_I = - \omega_I$ and $\varphi^*(\omega_I) =
\overline{\omega_I}$.
We define the real $2$-forms $Q$ and $R$ on $X$ to be $\omega_I = Q +
\sqrt{-1}R$.
Then we have
$$\tau^* Q = -Q,\ \ \tau^* R = -R,\ \
\varphi^* Q = Q \ \ \text{and} \ \ \varphi^* R = -R.$$ Moreover,
we may assume $2(P\wedge P) = \omega_I \wedge
\overline{\omega_I}$, equivalently, $Q \wedge Q = P \wedge P$.

We define the new complex structures $J$ and $K$ on $X$ by
$$
Q=g(J(\ ),\ )\ \ \text{and}\ \ R=g(K(\ ),\ ).
$$
(Then, $g$ is a hyperK\"ahler metric on $X$. See \cite{Huybrechts02},
Remark 2.4.)
By the above we can verify that
$\tau$ is anti-holomorphic and $\varphi$ is holomorphic for $(X,J)$.
We see $\omega_J := R + \sqrt{-1}P$ is
a nowhere vanishing holomorphic $2$-form on $(X,J)$.
We have $\varphi^*(\omega_J)= - \omega_J$, i.e.,
$\varphi$ is non-symplectic on the K3 surface $(X,J)$.
We get a new real K3 surface with a non-symplectic involution
$((X,J),\varphi,\tau)$.

Then the fixed point set $X_\varphi(\br)$ of $\varphi$
is a complex 1-dimensional submanifold of $(X,J)$ (\cite{Nikulin81}, p.1424).
The fixed point set of $\tau$ in $X_\varphi(\br)$ is $A(\br)$.
Thus,
we can apply Theorem \ref{realcurve} to
the K3 surface $(X,J)$, the anti-holomorphic involution $\tau$ and
the complex curve $C:=X_\varphi(\br)$.
Thus we conclude that the following two conditions are equivalent.\\
\indent
(1)\ $cl(X_\varphi(\br))\cdot x \equiv 0 \MOD 2 \ \forall x \in L_\tau$.\\
\indent (2)\ $cl(A(\br))=0\ \in H_1(A;\bz/2\bz)$.

The condition (1) is equivalent to the condition $\delta_{\varphi S}=0$.
This completes the theorem.
\end{proof}

In our case, $A=A_0\amalg A_1$ where $A_0\cong\bp^1$ and $A_1$ has
genus 9. Thus, $A_0$ is always dividing, and $A$ is dividing if
and only if $A_1$ is dividing which makes sense also for
$A=s+A_1\in |-2K_{\bff_3}|$. Thus, applying Theorem \ref{div}, we
get in our case that
\begin{equation}
A_1\ is\ dividing\ if\ and\ only\ if\ either\ (\delta_F=1,
\delta_\varphi =0)\ or\ \delta_F=\delta_{\varphi F}=0.
\label{dividingF3}
\end{equation}

Let us consider real $\bff_3$. Since $s^2=-3$ is odd, it follows
that $s(\br)^2\equiv 1\mod 2$ in $H_1(\bff_3(\br);\bz/2\bz)$. It
follows that $\bff_3(\br)$ is non-empty, and it is a
non-orientable surface. Rational pencil of $\bff_3$ gives then a
pencil on $\bff_3(\br)$ by circles over the circle $s(\br)$. Thus,
$\bff_3(\br)$ is Klein bottle. Further we denote by $c$ a real (i.
e. over a point of $s(\br)$) fibre of the pencil. Then $c(\br)$ is
a circle. The surface $Y$ (the blow-up of $\bff_3$ at $A_1\cap s$)
has $Y(\br)$ which is a non-orientable surface of genus $2$ ($\chi
(Y(\br))=-1$).

We consider isotopy type of the curve $A(\br)=s(\br)\cup A_1(\br)$
in $\bff_3(\br)$. By an {\it oval}, we mean a connected component of
$A_1(\br)$ which bounds a disk in $\bff_3(\br)$. Each otherwise
connected component is called a {\it non-oval}.

\medskip

We know that $A_1\cdot s=1$ and $A_1\cdot c=3$. Moreover, we know
the geometric interpretation of invariants \eqref{geninvF3} in
terms of the double covering K3 from Sect. \ref{geomK3F3}. Using
this information, by elementary 2-dimensional topology
considerations, we obtained all possible pictures.

Let $o_1$ be the connected component (a circle) of $A_1(\br)$ which contains
the unique point $A_1\cap s$. The $o_1$ is a non-oval, and there are
two cases.

\begin{figure}
\begin{center}
\includegraphics[width=10cm]{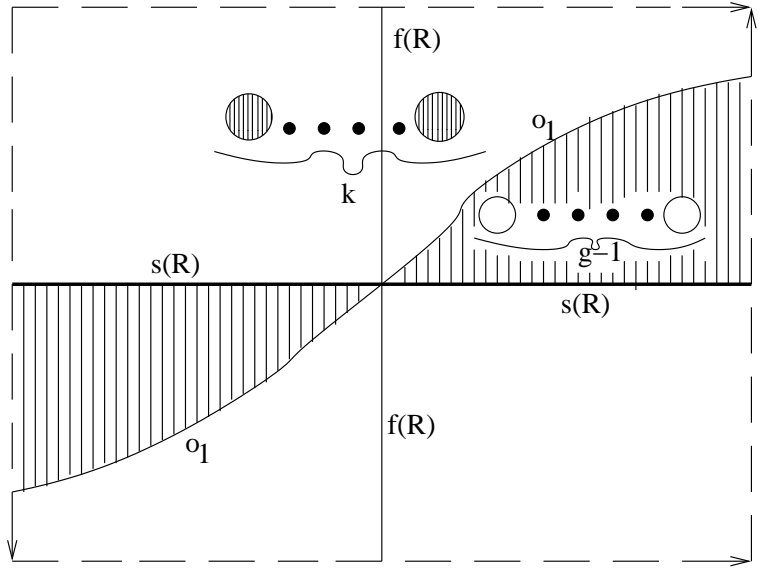}
\end{center}
\caption{The positive curve $A^+$ for $\delta_F=0$ and
$(r,a,\delta_{\varphi F})\not=(9,9,0)$ (see Figure \ref{diag0F3}),
here $g=(22-r-a)/2$, $k=(r-a)/2$.}
\label{Fig1F3}
\end{figure}

{\bf Case 1.}\ \ $A_1(\br)$ consists of the non-oval $o_1$ and
$g-1+k$ ovals. In case of $\delta_F=0$, the positive curve $A^+$
is shown as shadow in Figure \ref{Fig1F3}. As usual,
$g=(22-r-a)/2$ and $k=(r-a)/2$. This case does not include
invariants $(r,a,\delta_{\varphi F})=(9,9,0)$ which will appear in
the case 2 below. The negative part
$A^-=\overline{\bff_3(\br)-A^+}$ gives then the positive curve for
the related connected component of moduli with $\delta_F=1$. By
\eqref{relatedF3}, one can write down the corresponding related
invariants \eqref{geninvF3} of this positive curve.

In this case, the blow up $\pi(X_\varphi(\br))$ of $A^+$ consists
of a closed cylinder with $g-1$ holes, and $k$ closed disks. Two
boundary circles of the cylinder are $s(\br)$ and $o_1$. The blow
up $\pi(X_{\widetilde{\varphi}}(\br))$ of $A^-$ consists of a real
projective plane with $k+2$ holes, and $g-1$ disks. The boundaries
of two these holes are $s(\br)$ and $o_1$. Identifying the
boundaries of the corresponding disks and holes of
$\pi(X_\varphi(\br))$ and $\pi(X_{\widetilde{\varphi}}(\br))$, one
obtains $Y(\br)$.

\begin{figure}
\begin{center}
\includegraphics[width=10cm]{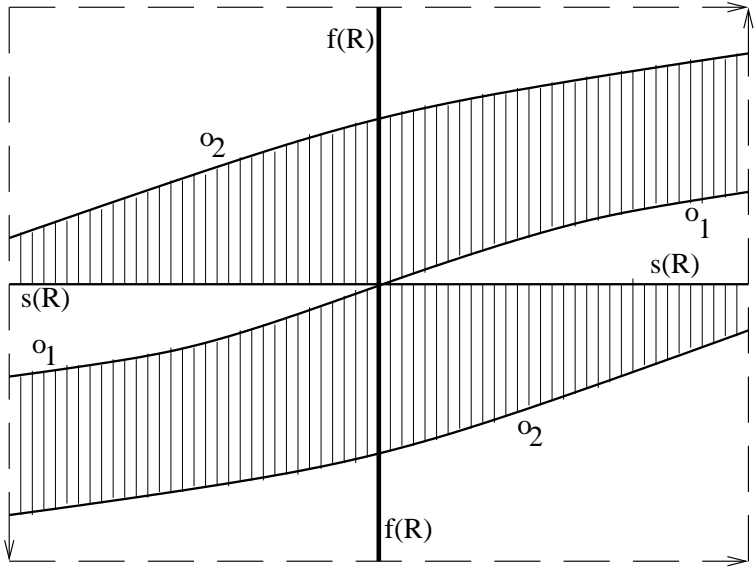}
\end{center}
\caption{The positive curve $A^+$ for $\delta_F=0$ and
$(r,a,\delta_{\varphi F})=(9,9,0)$ (see Figure \ref{diag0F3}).} \label{Fig2F3}
\end{figure}

{\bf Case 2.}\ \ $A_1(\br)$ consists of the non-oval $o_1$ and
another non-oval $o_2$. In case of $\delta_F=0$, the positive
curve $A^+$ is shown as shadow in Figure \ref{Fig2F3}. In this
case, $(r,a,\delta_{\varphi F})=(9,9,0)$. The negative part
$A^-=\overline{\bff_3(\br)-A^+}$ then corresponds to the positive
curve with the invariants $\delta_F=1$ and
$(r,a,\delta_{\varphi})=(10,8,0)$.

In this case, the blow up $\pi(X_\varphi(\br))$ of $A^+$ is a
closed disk with 2 holes. The boundary circles of the holes of the
disk are $s(\br)$ and $o_1$. The boundary of the disk is $o_2$.
The blow up $\pi(X_{\widetilde{\varphi}}(\br))$ of $A^-$ consists
of a closed cylinder and a M\"obius band. Two boundary circles of
the cylinder are $s(\br)$ and $o_1$. The boundary of the M\"obius
band is $o_2$.

These pictures and \eqref{dividingF3} show that the isotopy type
of $A^+$ together with dividedness determine the invariants
\eqref{geninvF3}. By Theorem \ref{themainF3} we finally get

\begin{theorem} The connected component of moduli of a positive
curve $A^+$ of a real non-singular curve
$A=s+A_1\in|-2K_{\bff_3}|$ is determined by the isotopy type of
$A^+$ and by dividedness by $A_1(\br)$ of $A_1(\bc)$.

All possibilities are presented in Figures \ref{diag0F3}
--- \ref{Fig2F3} and \eqref{dividingF3}. \label{connectedF3}
\end{theorem}

We also have similar result in terms of the blow up surface $Y$
(see its exact description above).

\begin{theorem} Let $Y$ be the blow up of $\bff_3$
at a point of the exceptional section of $\bff_3$.

The connected component of moduli of a positive curve $A^+$ of a
real non-singular curve $A\in|-2K_Y|$ is determined by the
topological  type of $A^+$ and by dividedness by $A(\br)$ of
$A(\bc)$.

All possibilities are presented in Figures \ref{diag0F3}
--- \ref{Fig2F3} and \eqref{dividingF3}. \label{connectedY}
\end{theorem}

One can consider other important applications of our results.

\begin{remark}
\label{remF4F3} Let us take the contraction of the exceptional
curve $F$ of $Y$. Then one obtains $\bff_4$. The image $A^\prime$
of $A$ gives then a curve $A^\prime \in |-2K_{\bff_4}|$ which is
non-singular except one quadratic (possibly degenerate) singular
point which is the image of $F$.

Thus, we obtain classification of connected components of moduli
(up to $Aut(\bff_4/\br)$) of real curves $A^\prime\in
|-2K_{\bff_4}|$ which are non-singular except one quadratic
singular point (the blow up at this point resolves this
singularity).
\end{remark}

We leave details to the interested reader.

\begin{remark}
\label{remZF3} Let us consider the contraction $\sigma:Y\to Z$ of
the curve $A_0\subset Y$ with $A_0^2=-4$. The surface $Z$ is a log
del Pezzo surface of index $\le 2$. Then image $\sigma(A_1)$ gives
a non-singular curve in $|-2K_Z|$, and vice versa. See
\cite{AlexeevNikulin88} for details. By classification of log del
Pezzo surfaces of index $\le 2$ in \cite{AlexeevNikulin88}, the
surface $Z$ is characterized by the properties: $Z$ has no
singular points except one log-terminal singular point of index
two having the type $K_1$ (it is the image of $A_0$), moreover the
genus $g(A_1)=9$.

Thus, we obtain classification of connected components of moduli
(up to $Aut(Z/\br)$) of real non-singular curves $A_1\in
|-2K_{Z}|$ where $Z$ is a real del Pezzo surface having no
singular points except one log-terminal singular point of index 2
of the type $K_1$, and the genus $g(A_1)=9$.
\end{remark}

We leave details to the interested reader.


\section{Connected components of moduli of real non-singular
curves in $|-2K_{\bff_2}|$} \label{secF2}

Here we consider the case of $\bff_2$. We shall see that this case
is more difficult, more interesting, very classical, and very
rich.

\subsection{Reduction to real K3 surfaces with
degenerate non-symplectic involution of type $(2,2,0)$}
\label{subsecK3F2}

Here we consider $Y=\bff_2$. Like above,  $c$ denotes the fibre of
the rational pencil of $\bff_2$, and $s$ its exceptional section.
We have $c^2=0$, $s^2=-2$ and $c\cdot s=1$. We have $-2K_{\bff_2}=
8c+4s$.

Let $A\in |-2K_{\bff_2}|$ be a non-singular curve. It is
irreducible of genus $9$. We have $A \cdot s=0$, and $A$ does not
intersect $s$.

Let $\pi:X\to Y=\bff_2$ be the double covering of $Y$ ramified in
$A\subset Y$. The $X$ is a K3 surface with the non-symplectic
involution $\tau$ of the double covering. Then $A=X^\tau$ is
identified with the fixed points set of the involution, and
$\bff_2=Y=X/\{1,\tau\}$ is the quotient space. We denote
$C=\pi^{-1}(c)$, $F=\pi^{-1}(s)$. Since $s\cong \bp^1$ and $A$
does not intersect $s$, the curve $F=F_1+F_2$ is disjoint union of
two non-singular rational curves which are conjugate by $\tau$. We
have $F_1^2=F_2^2=-2$, $F_1\cdot F_2=0$, $C\cdot F_1=C\cdot
F_2=1$. Here the numbering by $1$ and $2$ of the curves $F_1$ and
$F_2$ is not canonical, of course. The Gram matrix of $C$, $F_1$
and $F_2$ has the determinant $4$. Since there are no hyperbolic
even unimodular lattices of the rank $3$ (its signature must be
$0\mod 8$), it follows that $C$, $F_1$ and $F_2$ generate a
primitive sublattice in the Picard lattice $S_X$ of $X$.

Since $Y$ has Picard number $2$, it follows that $C$ and
$F=F_1+F_2$ generate the eigenvalue $1$ sublattice
$S=S_X^\tau=H_2(X(\bc);\bz)^\tau=\bz cl(C)+\bz cl(F)$ of $\tau$.
We have $r(S)=\rk S=2$, $S^\ast/S\cong(\bz/2\bz)^2$ is a
2-elementary group of the rank $a(S)=2$, we have $(x^\ast)^2\in
\bz$ for any $x^\ast\in S^\ast=S/2$. It follows that the invariant
$\delta(S)=0$. Thus, the main invariants of $\tau$ are $(r(S),
a(S),\delta(S))=(2,2,0)$. Moreover, in our case, all exceptional
curves on $X$ such that their image by the quotient  map $\pi:X\to
Y$ gives exceptional curves on $Y=\bff_2$, are the non-singular
rational curves $F_1$ and $F_2$ on $X$ such that $F_1^2=F_2^2=-2$,
$F_1\cdot F_2=0$ and $\tau(F_1)=F_2$. Other speaking, all
pre-images by $\pi$ of exceptional curves on $Y$ consist of the
exceptional curves $F_1$ and $F_2$.

Vice versa, let us consider a K3 surface with non-symplectic
involution $(X,\tau)$ where $S=(S_X)^\tau$ is a 2-elementary
lattice with the same invariants $(r(S),a(S),\delta(S))=(2,2,0)$.
By simple direct considerations, or by using general
classification results of \cite{AlexeevNikulin88}, there are two
cases for $X$ and the quotient  map $\pi:X\to Y=X/\{1,\tau\}$. In
both cases $Y=X/\{1,\tau\}$ is a rational surface with Picard
number $r(S)=2$, and $A=X^\tau$ gives a non-singular irreducible
curve $A=\pi(A)\in |-2K_Y|$ of genus $9$. Thus,  $\pi$ is a double
covering ramified in a non-singular curve $A\in |-2K_Y|$. In first
case, $Y$ has no exceptional curves. In second case, $Y$ has
exactly one exceptional curve $s$ with $s^2=-2$. By classification
of surfaces (e. g. see \cite{Shafarevich65}), then $Y\cong
\bff_0=\bp^1\times \bp^1$ in first case, and $Y=\bff_2$ in second
case. First case is called {\it non-degenerate.} Second case is
called {\it degenerate.}

The first, non-degenerate case, had been studied over $\br$ in
\cite{NikulinSaito05}.

The second, degenerate case, we consider here. We see that
non-singular curves $A\in |-2K_{\bff_2}|$ and K3 surfaces with
degenerate non-symplectic involution $(X,\tau)$ of type
$(2,2,0)$ are in one-to-one correspondence. They have the same
moduli and their connected components.

Now let us assume that $Y=\bff_2$ and a non-singular curve $A\in
|-2K_Y|$ are real. Let $\theta$ be the corresponding
anti-holomorphic involution on $Y$. Since $Y$ has the unique
rational pencil, $|c|$, and the unique exceptional curve, $s$,
then $\theta(cl(c))=-cl(c)$, $\theta(cl(s))=-cl(s)$, i. e.
$\theta$ acts as $(-1)$ in $H_2(Y(\bc);\bz)$. Let $\varphi$ be the
lift of $\theta$ on $X$. The $\varphi$ is defined uniquely up to
replacement by the related $\widetilde{\varphi}=\tau\varphi$.
Replacing $\varphi$ by $\widetilde{\varphi}$ if necessary, we can
always assume that $\varphi(F_1)=F_1$ and $\varphi(F_2)=F_2$. It
follows that $\varphi$ acts as $(-1)$ in the part of
$H_2(X(\bc);\bz)$ generated by $C$, $F_1$ and $F_2$. Thus the
action of $\varphi$ in $H_2(X(\bc);\bc)$ has at least 3 linearly
independent eigenvectors with the eigenvalue $(-1)$. If $\varphi$
has order $4$, then $\varphi^2=\tau$, and $\varphi$ must have
twenty linearly independent eigenvectors with eigenvalues $\pm
\sqrt{-1}$ in the eigenvalue $-1$ part $H_2(X(\bc);\bc)_\tau$ of
the action of $\tau$. Thus, we obtain more than $\rk
H_2(X(\bc);\bz)$ linearly independent eigenvectors for the action
of $\varphi$ which is impossible. This shows that $\varphi$ has
order $2$. The involutions $\varphi$ and $\tau$ obviously commute.

Thus, there exists a unique lift $\varphi$ of $\theta$ to an
anti-holomorphic involution $\varphi$ of $X$ which is determined
by the condition $\varphi(F_1)=F_1$ and $\varphi(F_2)=F_2$. Then
$\varphi$ acts as $-1$ in the $\bz$-submodule $[C,F_1,F_2]$ of
$H_2(X(\bc);\bz)$ generated by $C$, $F_1$ and $F_2$.
For the related involution
$\widetilde{\varphi}=\tau\varphi$ we have
$\widetilde{\varphi}(F_1)=F_2$ and $\widetilde{\varphi}
(F_2)=F_1$, $\widetilde{\varphi}$ acts as $(-1)$ in
$H_2(X(\bc);\bz)^\tau$, and
$\widetilde{\varphi}(cl(F_1)-cl(F_2))=cl(F_1)-cl(F_2)$ where
$cl(F_1)-cl(F_2)\in H_2(X(\bc);\bz)_\tau$.

Thus, classifying connected components of moduli of real K3
surfaces with degenerate non-symplectic involution
$(X,\tau,\varphi)$ of type $(2,2,0)$, we actually solve more
delicate problem of classification of connected components of
moduli of positive real non-singular curves $A^+$ where $A\in
|-2K_{\bff_2}|$ is non-singular. Here $A^+$ is determined by a
choice of one of two related anti-holomorphic involutions
$\varphi$ or $\widetilde{\varphi}$ of the double covering $X$. If
$Y(\br)\not=\emptyset$, we can identify $A^+$ with the closure of
one of two parts of $Y(\br)-A(\br)$ which is
$A^+=\pi(X_\varphi(\br))$. Then $A^-=\overline{Y(\br)-A^+}$ is
$A^-=\pi(X_{\widetilde{\varphi}}(\br))$.

If $Y(\br)\not=\emptyset$, then $s(\br)\not=\emptyset$, and the
part $A^+$ corresponding to the choice of $\varphi$ above is
determined by the condition $s(\br)\subset A^+$. Further we always
assume that the positive curve $A^+$ is chosen in this way: it
corresponds to $\varphi$ which acts on $[C,F_1,F_2]$ as $(-1)$.

\subsection{Enumeration of connected components of moduli of
real K3 surfaces with degenerate non-symplectic involution of
type $(2,2,0)$} \label{subsecModinvF2}

We correspond to a real K3 surface with non-symplectic involution
$(X,\tau,\varphi)$ (or $(X,\tau,\wvarphi)$) of type $(2,2,0)$ the
isomorphism class of the action $(L,\tau,\varphi)$ of the
involutions $\tau$ and $\varphi$ on the homology lattice
$L=H_2(X(\bc);\bz)$ (which is an even unimodular lattice of
signature $(3,19)$) where $\tau$ has type $(2,2,0)$ and $\varphi$
has hyperbolic $L^\varphi$. Here $\tau$ and $\varphi$ must act on
$[C,F_1,F_2]$ as has been described above: $\tau(cl(C))=cl(C)$,
$\tau(cl(F_1))=cl(F_2)$, $\tau(cl(F_2))=cl(F_1)$ and $\varphi$ is
$(-1)$ on $[C,F_1,F_2]$. Respectively $\wvarphi=\tau\varphi$ acts
on $[C,F_1,F_2]$ as follows: $\wvarphi(cl(C))=-cl(C)$,
$\wvarphi(cl(F_1))=-cl(F_2)$, $\wvarphi(cl(F_2))=-cl(F_1)$. Any
action with these properties takes place for some
$(X,\tau,\varphi)$ (or $(X,\tau,\wvarphi)$).

Two such actions $(L,\tau,\varphi)$ and
$(L^\prime,\tau^\prime,\varphi^\prime)$ are called isomorphic if
there exists an isomorphism $f:L\to L^\prime$ of lattices such
that $f\tau=\tau^\prime f$, $f\varphi=\varphi^\prime f$,
$f(cl(C))=cl(C)$, $f(cl(F_1))=cl(F_1)$ and $f(cl(F_2))=cl(F_2)$
(i. e. $f$ is identical on $[C,F_1,F_2]$. Replacing $f$ by
$f\tau$, we can always more generally assume that
$f(\{cl(F_1),cl(F_2)\}=\{cl(F_1),cl(F_2)\}$.

Since our case is degenerate (because $Y$ has real exceptional
curves with square $(-2)$), we cannot apply results of
(\cite{NikulinSaito05}, Sect. 2.1) where only non-degenerate case
had been considered. In spite of that, in our very special case
similar result is valid.

\begin{theorem} The connected component of moduli of a
real K3 surface $(X,\tau,\varphi)$ (or
$(X,\tau,\wvarphi=\tau{\varphi})$ ) with degenerate non-symplectic involution
of type $(2,2,0)$ is determined by the isomorphism class of
the action of $\tau$ and $\varphi$ (or $\wvarphi$) on the homology
lattice $L=H_2(X;\bz)$. \label{theactionF2}
\end{theorem}

\Proof We have $L^\tau=\bz cl(C)+\bz (cl(F_1)+cl(F_2))$. Its
orthogonal complement in $[C,F_1,F_2]$ is $\bz f$ where
$f=cl(F_1)-cl(F_2)$ and $f^2=-4$. We denote by $T$ the orthogonal
complement to $[C,F_1,F_2]$ in $L$.

Assume that $\delta\in L$ satisfies the properties: $\delta\perp
L^\tau$, the sublattice in $L$ generated by $[C,F_1,F_2]$ and
$\delta$ is hyperbolic, and $\delta^2=-2$.

Then $\delta=\alpha f+\beta \delta_2$ where $\alpha$, $\beta$ are
rational numbers, $\delta_2\in T$ is primitive and $\delta_2^2<0$.

We claim that then either $\alpha=0$ and $\delta=\pm \delta_2\in
T$, or $\alpha=\pm (1/2)$, $\beta=\pm (1/2)$ and
$f^2=\delta_2^2=-4$.

Really, if $\alpha=0$, then $\delta=\beta\delta_2\in T$ is
primitive, and then  $\beta=\pm 1$.

Assume that $\alpha\not=0$. We have $F_1\cdot \delta=-2\alpha\in
\bz$, and $\alpha=a/2$ where $a\in \bz$ and $a\not=0$. Then
$\delta=(a/2)f+\beta\delta_2$ and $2\delta-a f=2\beta\delta_2\in
L$. It follows that $2\beta\in \bz$ and $\beta=b/2$ where $b\in
\bz$. Thus $\delta=(a/2)f+(b/2)\delta_2$ where $a$ and $b$ are
non-zero integers. We have $\delta^2=-2=-a^2+(b^2/4) \delta_2^2$
where $\delta_2^2<0$ is even. It follows that $\delta_2^2=-4$ and
$a=\pm 1$, $b=\pm 1$.

Thus, $\delta=(\pm f\pm \delta_2)/2$ where $f\in [C,F_1,F_2]$ and
$\delta_2\in T$ are elements with square $-4$. We have
$[C,F_1,F_2]\cdot \delta=[C,F_1,F_2]\cdot f/2\in \bz$ and $T\cdot
\delta = T\cdot \delta_2/2\in \bz$. Thus, $[C,F_1,F_2]\cdot
f\equiv 0\mod 2$ and $T\cdot \delta_2\equiv 0\mod 2$. Thus,
$f^\prime=\pm f$ is a root of $[C,E_1,E_2]$ with square $(-4)$,
$\delta_2^\prime=\pm \delta_2$ is a root of $T$ with square
$(-4)$, and $\delta=(f^\prime+\delta_2^\prime)/2$ is the half sum
of these roots. Moreover, $\delta^\prime
=(f^\prime-\delta_2^\prime)/2\in L$ and $(\delta^\prime)^2=-2$.

If the mirror $\delta^\perp$ of $\delta$ divides the corresponding
period domain of $(X,\tau,\varphi)$ of our type in two disjoint
half-spaces, they will be conjugate by the reflection in $\delta$
in the first case, and by the composition of reflections in
$\delta$ and $\delta^\prime$ in the second case. It is easy to see
that these automorphisms will conjugate periods of isomorphic K3
surfaces with non-symplectic involution.

The remaining considerations are the same as for the
non-degenerate case in the proof of Theorem 1 in
\cite{NikulinSaito05}.

\QED

\begin{remark}
\label{remrootinvcompF2} Similar (properly formulated) result can
be proved for real K3 surfaces with degenerate non-symplectic
involution and the root invariant $(K,\xi)$ where $K$ is a
2-elementary root lattice (i. e. it is orthogonal sum of $A_1$,
$D_{2k}$, $E_7$, $E_8$). See \cite{AlexeevNikulin88} for
terminology and classification results. For our $\bff_2$ case, $K$
is $A_1$. For non-degenerate case, $K$ is zero.

Unfortunately, it seems, similar to Theorem \ref{theactionF2}
result holds only for very special classes of K3 surfaces with
degenerate non-symplectic involution. It would be interesting to
find all these nice cases.
\end{remark}

Thus, connected components of moduli of positive real non-singular
curves $A\in |-2K_Y|$ (further always $Y=\bff_2$) and,
equivalently, of real K3 surfaces with degenerate non-symplectic
involution $(X,\tau,\varphi)$ of type $(2,2,0)$ with our choice of
$\varphi$ between two related anti-holomorphic involutions are in
one-to-one correspondence with isomorphism classes of triplets
\begin{equation}
(L,\tau,\varphi) \label{triplinvF2}
\end{equation}
described above.

They can be enumerated using general results of \cite{Nikulin83}
about classification of integral involutions with conditions. On
the other hand, this case (like $\bff_3$ above) can be also
reduced to classification of integral polarized involutions (see
Sect. 3 in \cite{Nikulin79}) which is much simpler since all
necessary calculations are already done in (\cite{Nikulin79},
Theorem 3.4.3). We prefer the second possibility.

We can decompose $[C,F_1,F_2]=[C,F_1]\oplus [\widetilde{F}_1]$
where $\widetilde{F}_1=2C+(F_1-F_2)$. Here $[C,F_1]=U$ is even
unimodular lattice of signature $(1,1)$ and
$[\widetilde{F}_1]=\langle -4 \rangle$ since $(2C+F_1-F_2)^2=-4$

A triplet $(L,\tau,\varphi)$ from \eqref{triplinvF2} defines an
integral polarized involution
\begin{equation}
\label{intpolinvF2}
 (L_1,\varphi_1,cl(\widetilde{F}_1)).
\end{equation}
Here $L_1$ is the orthogonal complement to $U$ in $L$ (thus, it is
an even unimodular lattice of signature $(2,18)$,
$\varphi_1=\varphi\,|\,L_1$ and $L_1^{\varphi_1}$ is hyperbolic,
moreover, $\varphi_1(cl(\widetilde{F}_1))=-cl(\widetilde{F}_1)$ and
$\widetilde{F}_1^2=-4$.

Taking orthogonal sum of \eqref{intpolinvF2} with $U=[C,F_1]$, one
obtains the triplet \eqref{triplinvF2} where $\tau$ acts as
$\tau(cl(C))=cl(C)$, $\tau(cl(F_1))=cl(F_2)$, $\tau(cl(F_2))=cl(F_1)$, and
$\tau$ is $(-1)$
on $\widetilde{F_1}^\perp$ in $L_1$. Thus, their isomorphism
classes are in one to one correspondence.


Genus invariants of integral polarized involutions
\eqref{intpolinvF2} (and then \eqref{triplinvF2}) were completely
classified in Theorem 3.4.3 of \cite{Nikulin79}. Below we present
these results. (One should multiply by $(-1)$ the pairing in $L_1$
and apply Theorem 3.4.3 of \cite{Nikulin79} to $l_{(+)}=18$,
$l_{(-)}=2$, $t_{(-)}=1$, $t_{(+)}=r-1$ and $n=4$.)

The complete genus invariants of \eqref{intpolinvF2} and then
\eqref{triplinvF2} are
\begin{equation}
(r,a,\delta_\varphi;\delta_F,\delta_{\varphi F}).
 \label{geninvF2}
\end{equation}
Here $r=\rk L^\varphi\in \bn$;
$((L^\varphi)^\ast/L^\varphi)=(\bz/2\bz)^a$ where $a\ge 0$ is an
integer; $\delta_\varphi \in \{0,1\}$ is equal to $0$ if and only
if $x\cdot \varphi(x)\equiv 0\mod 2$ for any $x\in L$. They are
all invariants of the corresponding pair $(L,\varphi)$. Here
$\delta_{F}\in \{0,1\}$ is equal to $0$ if and only if $F \cdot
L_\varphi\equiv 0\mod 2$; here $\delta_{\varphi F}\in \{0,1\}$ is
equal to $0$ if and only if $x\cdot \varphi(x)\equiv x\cdot F$ for
any $x\in L$. Here and in what follows, we always denote by
$L^\varphi$ and $L_\varphi$ the eigenvalue $1$ and $-1$ parts
respectively of the action of an involution $\varphi$ on a module
$L$.

Strictly speaking, for invariants of integral polarized
involutions \eqref{intpolinvF2}, one should use $\widetilde{F}_1$
instead of $F$, $\varphi_1$ instead of $\varphi$,  and $L_1$
instead of $L$. But, $\widetilde{F}_1\equiv F=F_1+F_2\mod 2$ in
$L$, and the invariants are the same:
$\delta_F=\delta_{\widetilde{F}_1}$,
$\delta_\varphi=\delta_{{\varphi}_1}$, and $\delta_{\varphi
F}=\delta_{\varphi_1 \widetilde{F}_1}$.

The invariants \eqref{geninvF2} must satisfy the following
relations which are sufficient and necessary for existence of a
triplet \eqref{intpolinvF2} or the corresponding triplet
\eqref{triplinvF2} with these invariants.

\noindent {\bf 0. Conditions on $(r,a,\delta_\varphi )$:}

\noindent (1) $1\le r\le 18$, $0\le a\le \min \{r,\ 20-r\}$;

\noindent (2) $r+a\equiv 0\mod 2$; if $\delta_\varphi=0$, then
$r\equiv 2\mod 4$;

\noindent (3) if $a=0$, then $(\delta_\varphi=0,\ r\equiv 2\mod
8)$;

\noindent (4) if $a=1$, then $r\equiv 1,\,3\mod 8$;

\noindent (5) if $(a=2,\  r\equiv 6\mod 8)$, then
$\delta_\varphi=0$;

\noindent (6) if $(a=r,\ \delta_\varphi=0)$, then $r\equiv 2\mod
8$;

\noindent (7) if $(a=20-r, \delta_\varphi=0)$, then $r\equiv 2\mod
8$.

\medskip

\noindent{\bf I. Conditions on $\delta_F$ and $\delta_{\varphi
F}$:}

\noindent {\bf General conditions:}

\noindent (1) if $\delta_{\varphi F}=0$, then $(\delta_F=0,\
\delta_{\varphi }=1,\ r\equiv 0 \mod 4)$.

\medskip

\noindent {\bf Relations near the boundary $a=20-r$:}

\noindent (2) if $a=20-r$, then $\delta_F=0$;

\noindent (3) if $(a=20-r,\ \delta_{\varphi F}=0)$, then $r\equiv
0\mod 8$;

\noindent (4) if $(a=18-r,\ \delta_F=1,\ \delta_\varphi=0)$, then
$r\equiv 2\mod 8$.

\medskip

\noindent {\bf Relations near the boundary $a=0$:}

\noindent (5) if $a=0$, then $\delta_F=1$;

\noindent (6) if $a=1$, then $\delta_F=1$;

\noindent (7) if $(a=2,\ \delta_F=0,\ r\equiv 2\mod 8)$, then
$\delta_{\varphi}=0$;

\noindent (8) if $(a=2,\ \delta_F=0,\ r\equiv 0\mod 4)$, then
$\delta_{\varphi F}=0$.

\medskip

\begin{figure}
\begin{picture}(200,140)
\put(66,118){\circle{7}} \put(71,116){{\tiny means
$\delta_\varphi=0$ and $\delta_{\varphi F}=1$},}
\put(66,110){\circle{5}} \put(71,108){{\tiny means
$\delta_{\varphi}=1$ and $\delta_{\varphi F}=0$,}}
\put(66,102){\circle*{3}} \put(71,100){{\tiny means
$\delta_\varphi=\delta_{\varphi F}=1$}}

\multiput(8,0)(8,0){20}{\line(0,1){94}}
\multiput(0,8)(0,8){11}{\line(1,0){170}}
\put(0,0){\vector(0,1){100}} \put(0,0){\vector(1,0){180}} \put(
6,-10){{\tiny $1$}} \put( 14,-10){{\tiny $2$}} \put(
22,-10){{\tiny $3$}} \put( 30,-10){{\tiny $4$}} \put(
38,-10){{\tiny $5$}} \put( 46,-10){{\tiny $6$}} \put(
54,-10){{\tiny $7$}} \put( 62,-10){{\tiny $8$}} \put(
70,-10){{\tiny $9$}} \put( 76,-10){{\tiny $10$}} \put(
84,-10){{\tiny $11$}} \put( 92,-10){{\tiny $12$}}
\put(100,-10){{\tiny $13$}} \put(108,-10){{\tiny $14$}}
\put(116,-10){{\tiny $15$}} \put(124,-10){{\tiny $16$}}
\put(132,-10){{\tiny $17$}} \put(140,-10){{\tiny $18$}}
\put(148,-10){{\tiny $19$}} \put(156,-10){{\tiny $20$}}

\put(-8, -1){{\tiny $0$}} \put(-8,  7){{\tiny $1$}} \put(-8,
15){{\tiny $2$}} \put(-8, 23){{\tiny $3$}} \put(-8, 31){{\tiny
$4$}} \put(-8, 39){{\tiny $5$}} \put(-8, 47){{\tiny $6$}} \put(-8,
55){{\tiny $7$}} \put(-8, 63){{\tiny $8$}} \put(-8, 71){{\tiny
$9$}} \put(-10, 79){{\tiny $10$}} \put(-10, 87){{\tiny $11$}}
\put( -2,114){{\footnotesize $a$}} 
\put(186, -2){{\footnotesize $r$}} 

\put( 16,16){\circle{7}} \put( 48,16){\circle{7}} \put(
80,16){\circle{7}} \put(112,16){\circle{7}}
\put(144,16){\circle{7}} \put( 48,32){\circle{7}} \put(
80,32){\circle{7}} \put(112,32){\circle{7}} \put(
80,48){\circle{7}}
\put(80,64){\circle{7}}
\put( 80,80){\circle{7}}    
\put( 32,16){\circle{5}} \put( 32,32){\circle{5}}
\put( 64,16){\circle{5}} \put( 64,32){\circle{5}} \put(
64,48){\circle{5}} \put( 64,64){\circle{5}}
\put( 96,16){\circle{5}} \put( 96,32){\circle{5}} \put(
96,48){\circle{5}}
\put(128,16){\circle{5}} \put(128,32){\circle{5}}
\put( 24,24){\circle*{3}} \put( 40,24){\circle*{3}} \put(
56,24){\circle*{3}} \put( 72,24){\circle*{3}} \put(
88,24){\circle*{3}} \put(104,24){\circle*{3}}
\put(120,24){\circle*{3}} \put(136,24){\circle*{3}} \put(
32,32){\circle*{3}} \put( 48,32){\circle*{3}} \put(
64,32){\circle*{3}} \put( 80,32){\circle*{3}} \put(
96,32){\circle*{3}} \put(112,32){\circle*{3}}
\put(128,32){\circle*{3}} \put( 40,40){\circle*{3}} \put(
56,40){\circle*{3}} \put( 72,40){\circle*{3}} \put(
88,40){\circle*{3}} \put(104,40){\circle*{3}}
\put(120,40){\circle*{3}} \put( 48,48){\circle*{3}} \put(
64,48){\circle*{3}} \put( 80,48){\circle*{3}} \put(
96,48){\circle*{3}} \put(112,48){\circle*{3}} \put(
56,56){\circle*{3}} \put( 72,56){\circle*{3}} \put(
88,56){\circle*{3}} \put(104,56){\circle*{3}} \put(
64,64){\circle*{3}} \put( 80,64){\circle*{3}} \put(
96,64){\circle*{3}} \put( 72,72){\circle*{3}} \put(
88,72){\circle*{3}}\put( 80,80){\circle*{3}}   
\end{picture}
\caption{$\bff_2$: The complete list of the invariants
\ref{geninvF2} (i. e. in terms of $\varphi$) with $\delta_F=0$ of
the connected components of moduli.} \label{diag0F2}
\end{figure}
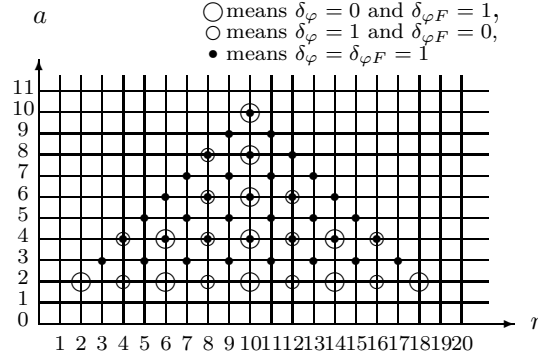

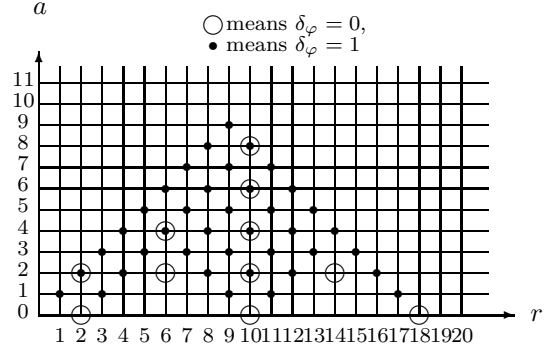
\begin{figure}
\begin{picture}(200,140)
\put(66,110){\circle{7}} \put(71,108){{\tiny means
$\delta_\varphi=0$,}} \put(66,102){\circle*{3}}
\put(71,100){{\tiny means $\delta_\varphi=1$}}

\multiput(8,0)(8,0){20}{\line(0,1){94}}
\multiput(0,8)(0,8){11}{\line(1,0){170}}
\put(0,0){\vector(0,1){100}} \put(0,0){\vector(1,0){180}} \put(
6,-10){{\tiny $1$}} \put( 14,-10){{\tiny $2$}} \put(
22,-10){{\tiny $3$}} \put( 30,-10){{\tiny $4$}} \put(
38,-10){{\tiny $5$}} \put( 46,-10){{\tiny $6$}} \put(
54,-10){{\tiny $7$}} \put( 62,-10){{\tiny $8$}} \put(
70,-10){{\tiny $9$}} \put( 76,-10){{\tiny $10$}} \put(
84,-10){{\tiny $11$}} \put( 92,-10){{\tiny $12$}}
\put(100,-10){{\tiny $13$}} \put(108,-10){{\tiny $14$}}
\put(116,-10){{\tiny $15$}} \put(124,-10){{\tiny $16$}}
\put(132,-10){{\tiny $17$}} \put(140,-10){{\tiny $18$}}
\put(148,-10){{\tiny $19$}} \put(156,-10){{\tiny $20$}}

\put(-8, -1){{\tiny $0$}} \put(-8,  7){{\tiny $1$}} \put(-8,
15){{\tiny $2$}} \put(-8, 23){{\tiny $3$}} \put(-8, 31){{\tiny
$4$}} \put(-8, 39){{\tiny $5$}} \put(-8, 47){{\tiny $6$}} \put(-8,
55){{\tiny $7$}} \put(-8, 63){{\tiny $8$}} \put(-8, 71){{\tiny
$9$}} \put(-10, 79){{\tiny $10$}} \put(-10, 87){{\tiny $11$}}
\put( -2,114){{\footnotesize $a$}} 
\put(186, -2){{\footnotesize $r$}} 

\put( 16, 0){\circle{7}} \put( 80, 0){\circle{7}} \put(144,
0){\circle{7}} \put( 16,16){\circle{7}} \put( 48,16){\circle{7}}
\put( 80,16){\circle{7}} \put(112,16){\circle{7}}
\put( 48,32){\circle{7}} \put(
80,32){\circle{7}}
\put(80,48){\circle{7}}
\put( 80,64){\circle{7}}

\put(  8, 8){\circle*{3}} \put( 24, 8){\circle*{3}} \put( 72,
8){\circle*{3}} \put( 88, 8){\circle*{3}} \put(136,
8){\circle*{3}}
\put(16,16){\circle*{3}}
\put( 32,16){\circle*{3}} \put( 64,16){\circle*{3}} \put(
80,16){\circle*{3}} \put( 96,16){\circle*{3}}
\put(128,16){\circle*{3}}
\put(24,24){\circle*{3}} \put( 40,24){\circle*{3}} \put(
56,24){\circle*{3}} \put( 72,24){\circle*{3}} \put(
88,24){\circle*{3}} \put(104,24){\circle*{3}}
\put(120,24){\circle*{3}}
\put(32,32){\circle*{3}} \put( 48,32){\circle*{3}} \put(
64,32){\circle*{3}} \put( 80,32){\circle*{3}} \put(
96,32){\circle*{3}} \put(112,32){\circle*{3}}
\put( 40,40){\circle*{3}} \put(
56,40){\circle*{3}} \put( 72,40){\circle*{3}} \put(
88,40){\circle*{3}} \put(104,40){\circle*{3}}
\put( 48,48){\circle*{3}} \put(
64,48){\circle*{3}} \put( 80,48){\circle*{3}} \put(
96,48){\circle*{3}}
\put(56,56){\circle*{3}} \put( 72,56){\circle*{3}} \put(
88,56){\circle*{3}}
\put( 64,64){\circle*{3}} \put( 80,64){\circle*{3}}
\put( 72,72){\circle*{3}}
\end{picture}

\caption{$\bff_2$: The complete list of the invariants
\ref{geninvF2} (i. e. in terms of $\varphi$) with $\delta_F=1$
(then $\delta_{\varphi F}=1$) of the connected components of
moduli.} \label{diag1F2}
\end{figure}

It is easy to enumerate all invariants \eqref{geninvF2} satisfying
these conditions. They are all presented in Figure \ref{diag0F2}
for $\delta_F=0$, and in Figure \ref{diag1F2} for $\delta_F=1$.


Additional considerations show that the genus invariants
\eqref{geninvF2} define the isomorphism class.

\begin{lemma} The genus invariants \eqref{geninvF2} determine the
isomorphism class of \eqref{intpolinvF2} and of
\eqref{triplinvF2}. Thus, they determine the connected component
of moduli. \label{thegenisomclassF2}
\end{lemma}

\Proof If $20-r\ge 4+\delta_F$, the statement follows from Theorem
3.3.1 in \cite{Nikulin79}.

Let us assume that $20-r<4+\delta_F$. Equivalently, $r\ge
17-\delta_F$. Let us denote $h=cl(\widetilde{F}_1)$. We have
$h^2=-4$.

In all cases it is enough to prove (see \cite{Nikulin79} or Remark
1.6.2 in \cite{Nikulin83}) that the lattice
$K=(L_1)_{\varphi_1,F}$ which is the orthogonal complement to $h$
in $(L_1)_{\varphi_1}$ has the properties:
\begin{equation}
\begin{array}{l}
K\ \text{is unique in its genus and the canonical homomorphism}\\
O(K)\to O(q_K)\text{ is epimorphic. }
\end{array}
\label{uniqF2}
\end{equation}
Here $q_K$ denotes the discriminant quadratic form of $K$.

Further we denote by $E_8$ and $E_7$ the negative definite
lattices defined by the Dynkin diagrams ${\mathbb E}_8$ and
${\mathbb E}_7$ respectively. For a lattice $M$ and a rational
$t$, we denote by $M(t)$ the lattice which is obtained by
multiplication by $t$ of the form of $M$.

Assume $\delta_F=1$. Then $r\ge 16$. By Figure \ref{diag1F2}, we
should check the genus invariants $(r,a,\delta_\varphi)=(18,2,0)$,
$(17,1,1)$ and $(16,2,1)$.

If $(r,a,\delta_\varphi)=(18,2,0)$, then $(L_1)_{\varphi_1}$ is an
even unimodular hyperbolic lattice of signature $(1,1)$, and $h\in
(L_1)_{\varphi_1}$ is a primitive element with $h^2=-4$. Then
$K=h^\perp$ in $(L_1)_{\varphi_1}$ is isomorphic to $\langle
4\rangle$, and \eqref{uniqF2} is obvious.

If $(r,a,\delta_\varphi )=(17,1,1)$ then $(L_1)^{\varphi_1}\cong
U\oplus E_8\oplus E_7$, and $K$ is orthogonal complement to
$(L_1)^{\varphi_1}\oplus \bz h$ in $L_1$. Then $K$ has the genus
of $\langle 4 \rangle \oplus \langle -2 \rangle$. Equivalently,
$K(1/2)$ has the genus of $\langle 2 \rangle \oplus \langle -1
\rangle$. Thus, $K(1/2)$ can be obtained as odd orthogonal
complement to an element $\delta\in M$ with $\delta^2=-2$ where
$M$ is the odd unimodular lattice $\langle 1 \rangle \oplus 2
\langle -1 \rangle$.

If $(r,a,\delta)=(16,2,1)$, then the same considerations show that
$K(1/2)$ can be obtained as odd orthogonal complement to an
element $\delta\in M$ with $\delta^2=-2$ where $M=\langle 1
\rangle \oplus 3\langle -1 \rangle$.

Assume $\delta_F=0$. Then $r\ge 17$, and the genus invariants of
Figure \ref{diag0F2} give $(r,a,\delta_\varphi)=(18,2,0)$,
$(17,3,1)$.

If $(r,a,\delta_\varphi)=(18,2,0)$, then the lattice
$(L_1)_{\varphi_1}\cong U(2)$, and $K$ is the orthogonal
complement to $h\in U(2)$ with $h^2=-4$. Then $K\cong \langle 4
\rangle$, and \eqref{uniqF2} is obvious.

If $(r,a,\delta_\varphi)=(17,3,1)$, then $(L_1)_{\varphi_1}\cong
\langle 2 \rangle\oplus 2\langle -2 \rangle$. It follows that
$K(1/2)$ is odd orthogonal complement to an element $h\in \langle
1 \rangle \oplus 2\langle -1 \rangle$ with $h^2=-2$.

Thus, for a lattice $M=\langle 1 \rangle\oplus t\langle -1
\rangle$ where $t=2$ or $t=3$ we have to check the following
statement: The odd orthogonal complements $K$ to elements
$\delta\in M$ with $\delta^2=-2$ are all isomorphic, and the
canonical homomorphism $O(K)\to O(q_{K(2)})$ is epimorphic.
Exactly this statement had been checked in the proof of Theorem 21
of \cite{NikulinSaito05}.

It follows the statement.  \QED

Thus, we finally enumerated connected components of moduli we are
looking for, and we obtained the main result.

\begin{theorem} There are exactly 108 connected components of moduli of
real non-singular curves $A\in |-2K_{\bff_2}|$ up to
$Aut(\bff_2/\br)$, equivalently of related pairs of positive real
non-singular curves $A\in |-2K_{\bff_2}|$, equivalently of related
pairs of connected components of moduli of real K3 surfaces with
degenerate non-symplectic involution of type $(2,2,0)$.

They are labelled by their invariants
\eqref{geninvF2} satisfying relations 0.1--7, I.1--8,
enumerated in Figures \ref{diag0F2} and \ref{diag1F2}.
\label{themainF2}
\end{theorem}

In \cite{NikulinSaito05} some other invariants of real K3 surfaces with
non-symplectic involution of type $(2,2,0)$ were used. For further
considerations and applications it is important to rewrite invariants
\eqref{geninvF2} in terms of that invariants.

For an involution $(L,\tau,\varphi)$ from \eqref{triplinvF2} and
$S=[C,F]=L^\tau$ we define the following invariants. First, like
$\delta_F$, for any $x_-\in S$ we can define $\delta_{x_-}\in
\{0,1\}$. Here $\delta_{x_-}=0$ if and only if $x_-\cdot
L_\varphi\equiv 0\mod 2$. For example, $\delta_C=1$ because
$C\cdot F_1=1$. Similarly $\delta_{C+F}=1$ because $(C+F)\cdot
F_1=-1$. Thus, only $\delta_F$ can take two values: $0$ or $1$.

Now we define a subgroup $H\subset S/2S$  by
the condition:
\begin{equation}
x_-\mod 2S\in H\ \text{iff\ }\delta_{x_-}=0.
\end{equation}
Thus, for our $(L,\tau,\varphi)$ we have
$H=\{0\}$ if $\delta_F=1$, and $H=\bz/ 2\bz cl(F)$ if $\delta_F=0$.

We say that $\delta_{\varphi S}=0$ if there exists $v\in S\mod 2S$
such that $x\cdot \varphi(x)\equiv x\cdot v\mod 2$ for any $x\in
L$. Otherwise, $\delta_{\varphi S}=1$. If $\delta_{\varphi S}=0$,
then the element $v\mod 2S$ is defined uniquely, and it is called
the {\it characteristic element of $\varphi$.} If $\delta_{\varphi
S}=1$, then the characteristic element $v$ of $\varphi$ is not
defined. It is easy to see that  if $\delta_{\varphi S}=0$ and $v$
is the characteristic element, then $\delta_v=0$.

In our case, we then have
$$
(\delta_{\varphi S}=0, v=0)\ \text{iff}\  \delta_\varphi=0;
$$
and
$$
(\delta_{\varphi S}=0, v=cl(F))\ \text{iff}\  \delta_{\varphi
F}=0;
$$
$$
\delta_{\varphi S}=1\ \text{iff}\ \delta_\varphi=\delta_{\varphi
F}=1.
$$
Thus, all invariants \eqref{geninvF2} are equivalent to
\begin{equation}
(r,a,H,\delta_{\varphi S}, v)
\label{invphiold}
\end{equation}
where $v$ is defined only for $\delta_{\varphi S}=0$.

The invariants \eqref{invphiold} are important because the same
invariants are defined for the related involution
$\widetilde{\varphi}=\tau\varphi$ by replacing everywhere
$\varphi$ by $\widetilde{\varphi}$. In \cite{NikulinSaito05} their
relation had been found. Let us denote by $a_H$ the rank of $H$
over $\bz/2\bz$. On $S/2S$ the non-degenerate pairing $b_S$ is
defined by the condition $b_S(x+2S,y+2S)=x\cdot y/2\mod 2$. It is
skew-symmetric: $b_S(x+2S,x+2S)=0\mod 2$. We have (see (4.2) and
(4.3) in \cite{NikulinSaito05}):
\begin{equation}
\begin{array}{ll}
&\varphi(cl(F_1))=-cl(F_1),\ \varphi(cl(F_2))=-cl(F_2);\
\widetilde{\varphi}(cl(F_1))=-cl(F_2),\
\widetilde{\varphi}(cl(F_2))=-cl(F_1);\\
&r(\varphi)+r(\widetilde{\varphi})=20;\
a_{H(\varphi)}+a_{H(\widetilde{\varphi})}=2;\
a(\varphi)-a_{H(\varphi)}=a(\widetilde{\varphi})-a_{H(\widetilde{\varphi})};\\
&H(\widetilde{\varphi})=H(\varphi)^\perp\ \text{with respect to }b_S;\
\delta_{\widetilde{\varphi} S}=\delta_{\varphi S};\
v(\widetilde{\varphi})=v(\varphi).
\end{array}
\label{relatedF2}
\end{equation}
Thus, $H(\wvarphi)=[F]=(\bz/2\bz)cl(F)$ if $\delta_F(\varphi)=0$,
and $H(\wvarphi)=[S]=S/2S$ if $\delta_F(\varphi)=1$. Applying
\eqref{relatedF2}, from invariants of Figures \ref{diag0F2},
\ref{diag1F2} we obtain all invariants of $\wvarphi$.

 Below we consider geometric interpretation of the invariants
\eqref{geninvF2} in terms of K3 and curves on $\bff_2$.

\subsection{Geometric interpretation of the invariants \eqref{geninvF2}
in terms of K3} \label{geomK3F2}

For the anti-holomorphic involutions $\varphi$ and
$\widetilde{\varphi}$ we have the same geometric interpretations
as in Sect. \ref{geomK3F3} for $\bff_3$. One should only use the
invariants $r(\varphi)$, $a(\varphi)$, $\delta_\varphi$,
$\delta_{x }(\varphi)$, $\delta_{x \varphi}$, $x \in S$, for
$\varphi$, and $r(\widetilde{\varphi})$, $a(\widetilde{\varphi})$,
$\delta_{\widetilde{\varphi}}$, $\delta_{x}(\widetilde{\varphi})$,
$\delta_{x \widetilde{\varphi}}$, for $\widetilde{\varphi}$. For
example, $\delta_{\widetilde{\varphi}C}=0$ if and only if
$X_{\widetilde{\varphi}}(\br)\sim C(\bc)\mod 2$ in
$H_2(X(\bc);\bz)$.

Using results above, like for $\bff_2$ we then obtain

\begin{theorem} The connected component of moduli of a
real K3 surface with degenerate
non-symplectic involution $(X,\tau,\varphi)$
(or $(X,\tau,\widetilde{\varphi})$)of type $(2,2,0)$,
equivalently of a positive real non-singular curve $A\in
|-2K_{\bff_2}|$, is determined by the geometric invariants:
topology of the real part $X_\varphi(\br)$, if $F_\varphi(\br)\sim
0\mod 2$ in $H_1(X_\varphi(\br); \bz)$, if $X_\varphi (\br)\sim
0\mod 2$ in $H_2(X(\bc);\bz)$, if $X_\varphi (\br)\sim F(\bc)\mod
2$ in $H_2(X(\bc);\bz)$.
\label{thegeomK31F2}
\end{theorem}

\begin{theorem} The connected component of moduli of a real K3 surface
with degenerate non-symplectic involution $(X,\tau,\varphi)$
(or $(X,\tau,\widetilde{\varphi})$)of type $(2,2,0)$,
equivalently of a positive real non-singular curve $A\in
|-2K_{\bff_2}|$, is determined by the geometric invariants:
topology of the real part $X_{\widetilde{\varphi}}(\br)$,
if $C_{\widetilde{\varphi}}(\br)\sim
0\mod 2$ in $H_1(X_{\widetilde{\varphi}}(\br); \bz)$,
 if
$X_{\widetilde{\varphi}} (\br)\sim
0\mod 2$ in $H_2(X(\bc);\bz)$, if $X_{\widetilde{\varphi}} (\br)\sim F(\bc)\mod
2$ in $H_2(X(\bc);\bz)$.
\label{thegeomK32F2}
\end{theorem}

\subsection{Geometric interpretation of the invariants \eqref{geninvF2}
in terms of $A\in |-2K_{\bff_2}|$} \label{geomF2F2} By Theorem
\ref{div}, {\it the real curve $A=X^\tau$ is dividing, i. e.
$cl(A(\br))=0$ in $H_1(A(\bc);\bz/2\bz)$, if and only if
$\delta_{\varphi S}=0$ (equivalently, $\delta_{\wvarphi S}=0$).}

By results above, $X_\varphi(\br)$ and
$X_{\widetilde{\varphi}}(\br)$ are both empty and then
$\bff_2(\br)$ is empty if and only if the invariants
\eqref{geninvF2} are: $\delta_F=0$ and
$(r,a,\delta_\varphi)=(10,10,0)$. This gives only one connected
component of moduli of non-singular curves $A\in |-2K_{\bff_2}|$,
and two related connected components of moduli of positive curves.

Further we assume that $\bff_2(\br)\not=\emptyset$, equivalently,
the invariants \eqref{geninvF2} are different from
$(\delta_F,r,a,\delta_{\varphi})=(0,10,10,0)$.  See Figures
\ref{diag0F2} and \ref{diag1F2}. Then $s(\br)\not=\emptyset$, and
the rational pencil of $\bff_2$ gives a pencil over the circle
$s(\br)$ by circles $c(\br)$. Further we always take the rational
fibre to be real, i. e. over $s(\br)$. Since $s^2=-2$ is even, it
follows that $\bff_2(\br)$ is a torus where $s(\br)$ and $c(\br)$
are generators of the torus (i. e. of $H_1(\bff_2(\br);\bz)$).

A component of $A(\br)$ is called {\it oval} if it bounds a disk
in $\bff_2(\br)$. Otherwise, it is called non-oval.

We always consider the positive curve $A^+=\pi(X_\varphi(\br))$,
i. e. when $s(\br)\subset A^+$. Since $A$ does not intersect $s$
and $c\cdot A=4$, by elementary 2-dimensional topology
considerations, we have the following and only the following
pictures where we use the geometric interpretation of invariants
\eqref{geninvF2} above in terms of the double covering K3.
Moreover, to draw correctly $s(\br)$ in Figures \ref{Fig01F2} and
\ref{Fig02F2}, we use ellipsoid deformation from Sect.
\ref{subsecelldefF2} below.

\begin{figure}
\begin{center}
\includegraphics[width=10cm]{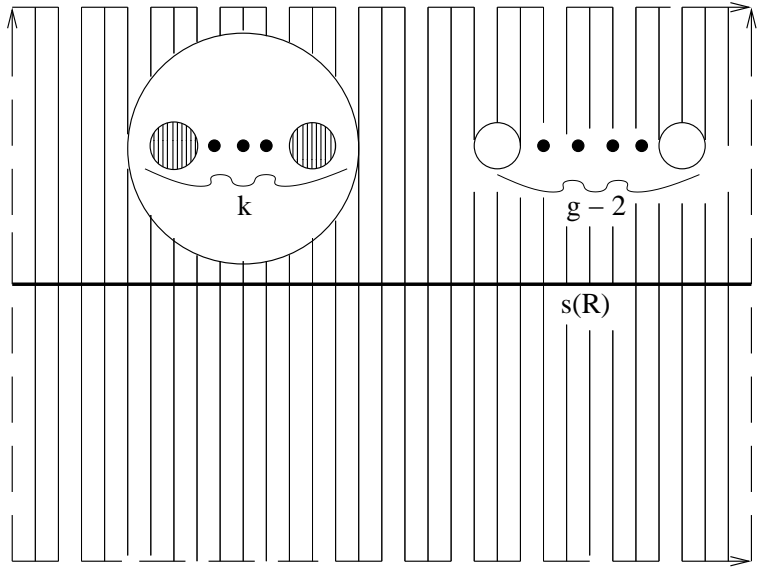}
\end{center}
\caption{$\bff_2$: The positive curve $A^+=\pi(X_\varphi(\br))$ for
$\delta_F=1$ and $(r,a,\delta_{\varphi})\not=(10,8,0)$, $(10,6,0)$
(see Figure \ref{diag1F2}),
here $g=(22-r-a)/2$, $k=(r-a)/2$.}
\label{Fig11F2}
\end{figure}

\begin{figure}
\begin{center}
\includegraphics[width=10cm]{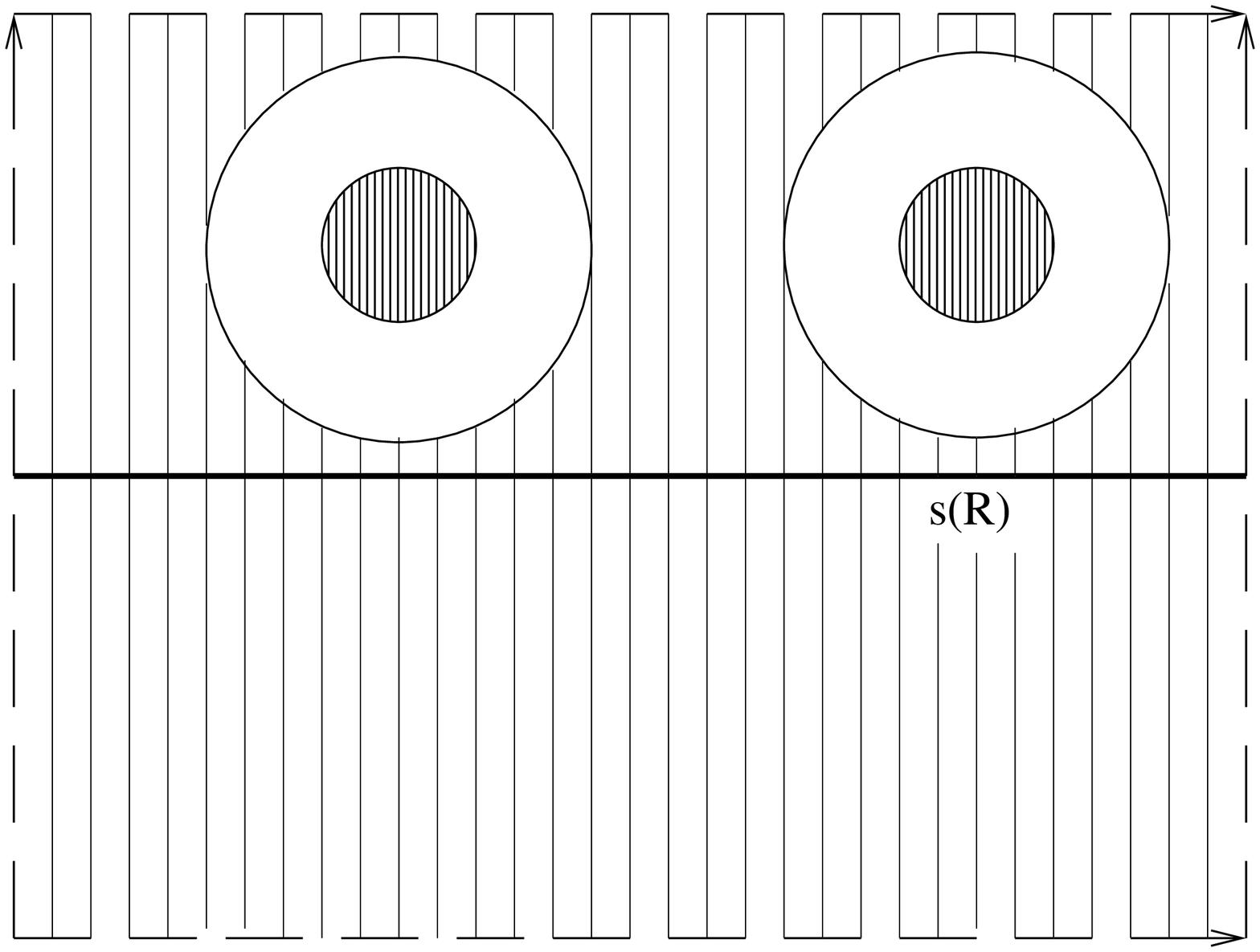}
\end{center}
\caption{$\bff_2$: The positive curve $A^+=\pi(X_\varphi(\br))$ for
$\delta_F=1$ and $(r,a,\delta_{\varphi})=(10,6,0)$
(see Figure \ref{diag1F2}).}
\label{Fig12F2}
\end{figure}

\begin{figure}
\begin{center}
\includegraphics[width=10cm]{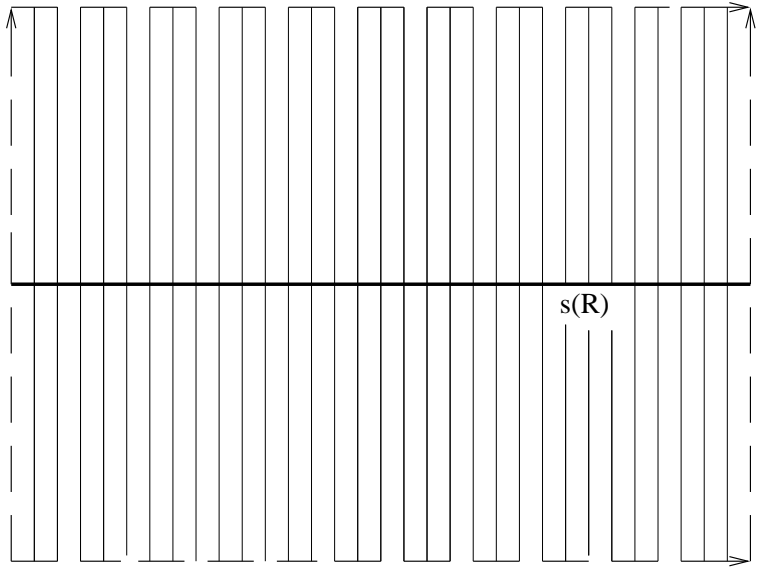}
\end{center}
\caption{$\bff_2$: The positive curve $A^+=\pi(X_\varphi(\br))$ for
$\delta_F=1$ and $(r,a,\delta_{\varphi})=(10,8,0)$ (see Figure \ref{diag1F2}).
Then $X_\varphi(\br)=(T_1)^2$ and $A(\br)=\emptyset$.}
\label{Fig13F2}
\end{figure}

{\bf Case 1.} $A(\br)$ has no non-ovals in
$\bff_2(\br)\not=\emptyset$. Then $\delta_F=1$ for the invariants
\eqref{geninvF2}, and the positive curve $A^+=\pi(X_\varphi(\br))$
is shown (as shadow) in Figures \ref{Fig11F2}, \ref{Fig12F2} and
\ref{Fig13F2}. All invariants \eqref{geninvF2} with $\delta_F=1$
are covered by this case.

{\bf Case 2.} $A(\br)$ has non-ovals. Then $\delta_F=0$ for the
invariants \eqref{geninvF2}, and the positive curve
$A^+=\pi(X_\varphi(\br))$ is shown in Figures \ref{Fig01F2},
\ref{Fig02F2} and \ref{Fig03F2}.
This case covers all invariants \eqref{geninvF2}
with $\delta_F=0$ except $(r,a,\delta_{\varphi})=(10,10,0)$ when
$\bff_2(\br)=\emptyset$. See Figure \ref{diag0F2} for the full
list of these invariants. If $\delta_F=0$ and the invariants
$(r,a,\delta_\varphi )\not=(10,10,0)$ nor $(10,8,0)$, then
$A(\br)$ has exactly two non-ovals, and $A^+$ is shown in Figures
\ref{Fig01F2}, \ref{Fig02F2}. If $(r,a,\delta_\varphi )=(10,8,0)$,
then $A(\br)$ consists of four non-ovals, and $A^+$ is shown in Figure
\ref{Fig03F2}.

Thus, we obtain

\begin{figure}
\begin{center}
\includegraphics[width=10cm]{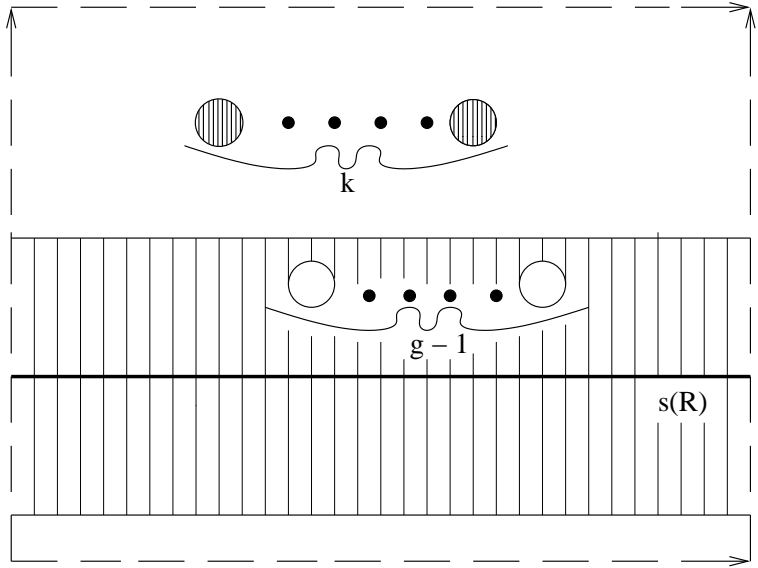}
\end{center}
\caption{$\bff_2$: The positive curve $A^+=\pi(X_\varphi(\br))$ for
$\delta_F=0$ and $(r,a,\delta_{\varphi})\not=(10,10,0)$, $(10,8,0)$,
and $(r,a,\delta_{\varphi F})\not=(8,8,0)$
(see Figure \ref{diag0F2}),
here $g=(22-r-a)/2$, $k=(r-a)/2$.}
\label{Fig01F2}
\end{figure}

\begin{figure}
\begin{center}
\includegraphics[width=10cm]{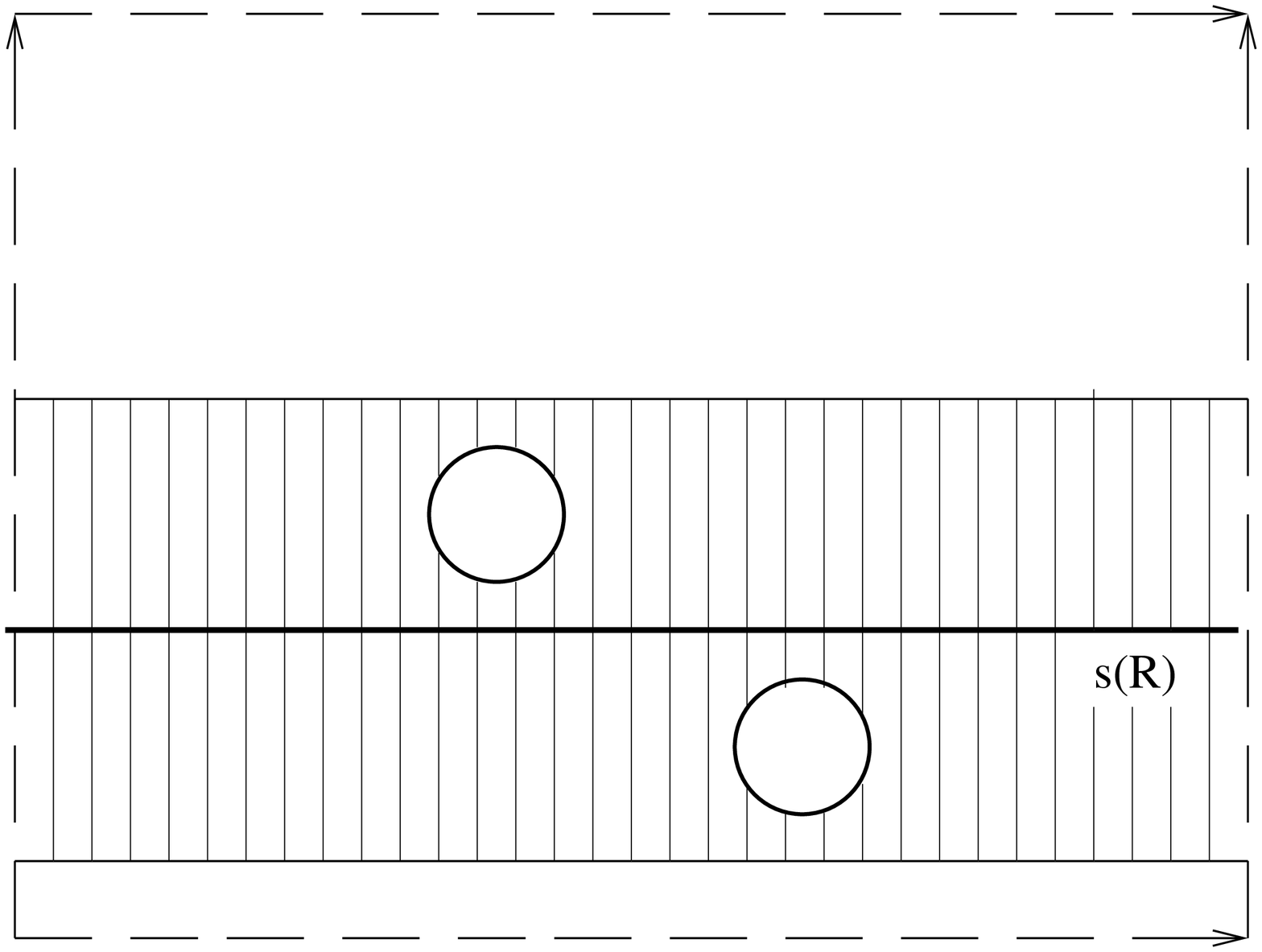}
\end{center}
\caption{$\bff_2$: The positive curve $A^+=\pi(X_\varphi(\br))$ for
$\delta_F=0$ and $(r,a,\delta_{\varphi F})=(8,8,0)$
(see Figure \ref{diag0F2}).}
\label{Fig02F2}
\end{figure}

\begin{figure}
\begin{center}
\includegraphics[width=10cm]{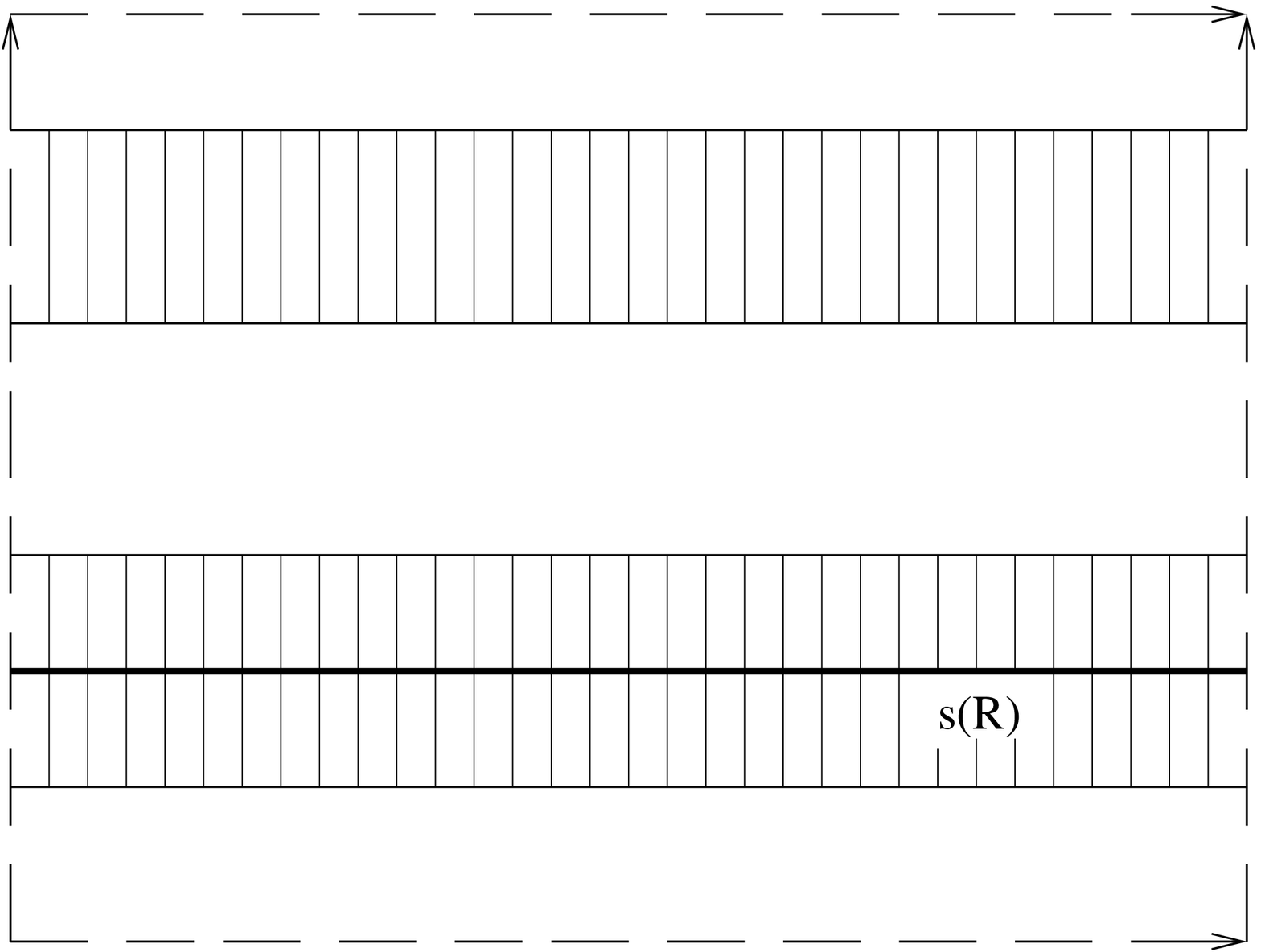}
\end{center}
\caption{$\bff_2$: The positive curve $A^+=\pi(X_\varphi(\br))$ for
$\delta_F=0$ and $(r,a,\delta_{\varphi})=(10,8,0)$
(see Figure \ref{diag0F2}).}
\label{Fig03F2}
\end{figure}

\begin{theorem} The connected component of moduli of a positive
curve $A^+$ of a real non-singular curve $A\in|-2K_{\bff_2}|$ is
determined by the following invariants: the isotopy type of $A^+$
in $\bff_2(\br)$, by dividedness by $A(\br)$ of $A(\bc)$, if
$s(\br)\subset A^+$. If $\bff_2(\br)=\emptyset$, there are exactly
two related connected components of moduli.

All possibilities are presented in Figures \ref{diag0F2}
---\ref{Fig03F2} and \ref{relatedF2}.
\label{theconnectedF2}
\end{theorem}

\newpage

\subsection{Application: Connected components of moduli of real non-singular
degree $8$ curves on a quadratic cone in $\bp^3$.}
\label{subsecapp1F2} We have $h_1=-K_{\bff_2}/2=2c+s$ has
$h_1^2=2$. The linear system $|h_1|$ gives the embedding of
$\bff_2$ as a quadratic cone in $\bp^3$ with the exceptional
section $s$ contracted to the vertex of the cone. The linear
system $|-2K_{\bff_2}|$ consists then of sections of the cone by
the degree $4$ forms, i. e. of the degree $8$ curves, which do not
pass through the vertex of the cone.

Let $(\xi_0:\xi_1:\xi_2:\xi_3)$ be usual real homogeneous
coordinates of $\bp^3$. The cone is given by the equation
$Q_2=\xi_1^2+\xi_2^2-\xi_3^2=0$ if $\bff_2(\br)\not=\emptyset$,
and by the equation $\xi_1^2+\xi_2^2+\xi_3^2=0$ if
$\bff_2(\br)=\emptyset$. The real curve $A$ is given by the
equation $F_4=0$ on the cone $Q_2=0$ which must be non-singular as
a complex algebraic curve. Here $F_4(\xi_0,\xi_1,\xi_2,\xi_3)$ is
a homogeneous degree four form with real coefficients which is
defined up to $\lambda(F_4+Q_2F_2)$ where $\lambda \in \br^\ast$
and $F_2$ a degree two form $F_2(\xi_0,\xi_1,\xi_2,\xi_3)$ with
real coefficients. The positive curve $A^+$ is given by the
inequality $F_4\ge 0$ (or $-F_4\ge 0$) on the cone such that $A^+$
contains the vertex $(1:0:0:0)$ of the cone.

Thus, we obtain classification of connected components of moduli
of real non-singular degree 8 curves on the real cone in $\bp^3$
up to the automorphism group of the cone which is the group
$PO(Q_2)$. It consists of projective real transformations which
multiply the form $Q_2$ by a constant.

We leave further geometric interpretation of our results above in
terms of degree 8 curves on a quadratic cone to the interested
reader.

\medskip

It is important and interesting to consider deformations of
non-singular degree 8 curves on a quadratic cone to non-singular
degree 8 curves on a non-degenerate quadric: One replaces the
degenerate quadratic form $Q_2$ defining the cone $Q_2=0$ (which
is $\bff_2$ with contracted $s$) by a closed non-degenerate
quadratic form $Q^\prime_2$, and the degree four form $F_4$
defining a non-singular curve $A$ on the cone ($Q_2=F_4=0$) by a
closed degree four form $F_4^\prime$. Then equations
$Q_2^\prime=F_4^\prime=0$ define a non-singular bidegree $(4,4)$
curve $A^\prime$ on the non-degenerate quadric $Q^\prime=0$ which
is $\bff_0$.

There are two very different types of these deformations. First
type when the quadric $Q_2^\prime=0$ is different from ellipsoid
(it is a hyperboloid or empty), e. g. the form
$Q_2=\xi_1^2+\xi_2^2-\xi_3^2$ is replaced by the form
$-\varepsilon \xi_0^2+\xi_1^2+\xi_2^2-\xi_3^2$ where
$\varepsilon>0$ is small. Second type when the quadric
$Q^\prime_2=0$ is ellipsoid, e. g. the form
$Q_2=\xi_1^2+\xi_2^2-\xi_3^2$ is replaced by the form
$Q_2^\prime=\varepsilon \xi_0^2+\xi_1^2+\xi_2^2-\xi_3^2$ where
$\varepsilon >0$ is small. Second type deformations, when the
quadric $Q_2^\prime=0$ is ellipsoid, is called {\it ellipsoid
deformation.} All other are called {\it hyperboloid deformations.}

Further we study these deformations in details.

\subsection{Application: Hyperboloid deformations of $A\in |-2K_{\bff_2}|$ to
$A\in |-2K_{\bff_0}|$.} \label{subsechypdefF2}

Let us consider a real K3 surface with degenerate
non-symplectic involution $(X,\tau,\varphi)$ (or $(X,\tau,\wvarphi)$)) of
type $(2,2,0)$. It defines then a connected component of moduli of
real K3 surfaces with non-symplectic involution of type
$(2,2,0)$ where we consider non-degenerate and degenerate
involutions all together. By \cite{NikulinSaito05}, the connected
component of moduli is also defined by the isomorphism class of
the action $(L,\tau,\varphi)$ with the condition $\varphi=-1$ on
$S=L^\tau=[C,F=F_1+F_2]$ for the hyperboloid case (for ellipsoid
case it is different, see below). We just forget about $F_1$ and
$F_2$, and are interested only on their sum
$cl(F)=cl(F_1)+cl(F_2)$. Moreover, connected components of moduli
of all involutions and non-degenerate ones coincide.

In \cite{NikulinSaito05} other generators for $S=L^\tau$ (related
with $\bff_0$) were used. They are $e_1=cl(C)$, $e_2=cl(C)+cl(F)$
(i. e. the effective isotropic generators of $S$ corresponding to
the structure $\bff_0=\bp^1\times \bp^1$) and
$h=e_1+e_2=2cl(C)+cl(F)$ (which is pre-image of hyperplane of
$\bp^3$). Thus $cl(F)\equiv e_1+e_2=h\mod 2$. Invariants of
connected components of moduli for hyperboloid case are the same
as we used  in Sect. \ref{subsecModinvF2}. See
\cite{NikulinSaito05}. Thus, hyperbolic deformations of
non-singular curves $A\in|-2K_{\bff_2}|$ to non-singular curves
$A\in |-2K_{\bff_0}|$ can be described symbolically as follows:
\begin{equation}
\begin{array}{ll}
(\bff_2;r(\varphi),a(\varphi),H(\varphi),\delta_{\varphi
S},v(\varphi))\Longrightarrow
(\bhh:r=r(\varphi),a=a(\varphi),H=H(\varphi),\delta_{\varphi
S},v=v(\varphi));\\
(\bff_2;r(\wvarphi),a(\wvarphi),H(\wvarphi),\delta_{\wvarphi
S},v(\wvarphi))\Longrightarrow
(\bhh:r=r(\wvarphi),a=a(\wvarphi),H=H(\wvarphi),\delta_{\varphi
S}=\delta_{\wvarphi S},v=v(\wvarphi)).
\end{array}
\label{defF2H2}
\end{equation}
where to the left we have invariants of a connected component of
moduli of positive real non-singular curves on $\bff_2$, and to
the right we have invariants of a connected component of moduli of
positive real non-singular curves on hyperboloid (see
\cite{NikulinSaito05}) obtained under hyperboloid deformation.

Using \eqref{defF2H2}, we can answer the following interesting
question. When there exist such deformation, and how many of them
do exist? Comparing invariants of connected components for
hyperboloid from \cite{NikulinSaito05} and for $\bff_2$ from Sect.
\ref{subsecModinvF2}, we obtain the result.

For $\bff_2$, the groups $H(\varphi)$ and $H(\wvarphi)$ cannot be
equal to $\bz/2\bz e_1$ or $\bz/2\bz e_2$. The group $H(\varphi)$
cannot be $S/2S$. The group $H(\wvarphi)$ cannot be zero. The
invariants $(r=12,a=8,H=[e_1+e_2],\delta_{\varphi S}=0,v=e_1+e_2)$
of hyperboloid cannot be obtained from $\varphi$ on $\bff_2$ (see
Figure \ref{diag0F2}) but they can be obtained from $\wvarphi$ on
$\bff_2$. Similarly,  the invariants $(r=8,a=8,H=[e_1+e_2],
\delta_{\varphi S}=0,v=e_1+e_2)$ of hyperboloid can be obtained
from $\varphi$ on $\bff_2$, but they cannot be obtained from
$\wvarphi$ on $\bff_2$. Using Sect. \ref{geomF2F2} and geometric
interpretation of the invariants for hyperboloid from
\cite{NikulinSaito05}, we can formulate the result in geometric
terms. Fortunately, here we don't have too many possibilities.

\begin{theorem}
\label{thedefF2hyperbF2}
If $\bff_0/\br$ is different from ellipsoid,
then a positive real non-singular curve $A\in |-2K_{\bff_0}|$
can be obtained by deformation of a positive real non-singular curve $A_0\in
|-2K_{\bff_2}|$ if and only if $\bff_0/\br$ is hyperboloid and
$A(\br)$ does not have connected components homological to any
line ($l_1$ or $l_2$ from different generators) of hyperboloid.
See Figures 10---13 in \cite{NikulinSaito05}.

This deformation is unique if either $A(\br)$ consists of only
ovals (including empty $A(\br)$) on hyperboloid (see Figures 10,
11 in \cite{NikulinSaito05}), or $A(\br)$ consists of two
non-ovals homological to $l_1+l_2$ and two ovals between these
non-ovals (then $A^+$ is connected) and $A$ is dividing. (See
Figure 12 in \cite{NikulinSaito05} for $g=3$, $k=0$ or $g=1$,
$k=2$.) In the last case, the connected component of the positive
curve of $A_0$ containing two ovals, must also contain $s(\br)$.

Otherwise, there are exactly two isomorphism classes of such deformations:
for one of them $s(\br)$ is contained in the positive curve of $A_0$, and for
another one $s(\br)$ is not contained in the positive curve of $A_0$.
\end{theorem}

\subsection{Application: Ellipsoid deformations of $A\in |-2K_{\bff_2}|$ to
$A\in |-2K_{\bff_0}|$.} \label{subsecelldefF2} We describe the
ellipsoid deformation of a positive real non-singular curve  $A\in
|-2K_{\bff_2}|$ in terms of periods of K3 surfaces using Global
Torelli Theorem for K3 surfaces \cite{PS71}.

For an $A\in |-2K_{\bff_2}|$, let $(X,\tau,\varphi)$ be the
corresponding real K3 surface with degenerate non-symplectic
involution of type $(2,2,0)$, and $(L,\tau,\varphi)$ its
action.

An element $\delta$ of a lattice $M$ with $\delta^2=-2$ defines
the reflection $s_\delta$ of the lattice $M$ by the formula $x\to
x+(x\cdot \delta)\delta$, $x\in M$. We have
$s_\delta(\delta)=-\delta$ and $\delta$ is identical in
$\delta^\perp$.

Since $F_1^2=F_2^2=-2$, they define reflections $s_{cl(F_1)}$ and
$s_{cl(F_2)}$ of $L$. Since $F_1\cdot F_2=0$, the reflections
$s_{cl(F_1)}$ and $s_{cl(F_2)}$ commute, and we obtain the
involution
$$
s_{[F_1,F_2]}=s_{cl(F_1)}s_{cl(F_2)}
$$
of $L$ with the eigenvalue $-1$ part $[F_1,F_2]$.

The involution $s_{[F_1,F_2]}$ obviously commutes with $\tau$ and
$\varphi$. Let us replace the involutions $\varphi$ and $\wvarphi$
on $L$ by the involutions (we apply twice the Lefschetz formula
for the vanishing cycle)
\begin{equation}
 \psi=s_{[F_1,F_2]}\varphi\,,\ \
 \wpsi=s_{[F_1,F_2]}\wvarphi=\tau\psi\ .
\label{psiF2}
\end{equation}
They coincide with $\varphi$, $\wvarphi$ in $[F_1,F_2]^\perp$
respectively. This geometrically means that $\psi$, $\varphi$ and
$\wpsi$, $\wvarphi$ are the same outside the exceptional curves
$F_1$ and $F_2$ which must take place for ellipsoid deformation.

The involution $\psi$ coincides with $\varphi$ in
$[C,F_1,F_2]^\perp$, and $\psi(cl(F_1))=cl(F_1)$,
$\psi(cl(F_2))=cl(F_2)$. As above, $h=2cl(C)+cl(F)$ generates the
orthogonal complement to $F$ in $S=[C,F]$, and $\psi(h)=-h$,
$\psi(cl(F))=cl(F)$. Thus, $\psi$ acts in $S=L^\tau$ as ellipsoid
involution (see \cite{NikulinSaito05}), and $(L,\tau,\psi)$
defines the connected component of moduli of curves $A\in
|-2K_{\bff_0}|$ where $\bff_0$ is ellipsoid, see
\cite{NikulinSaito05}.

For degenerate $(X,\tau,\varphi)$ its periods (e. g. see
\cite{Nikulin79} and Sect. 2.1 in \cite{NikulinSaito05}) are
equivalent to a pair of rays $(\br^+\omega_+,\br^+\omega_-)$ where
$\omega_+\in L^\varphi\otimes \br$ and $\omega_+^2>0$, and
$\omega_-\in L_{\varphi}\otimes \br$ where $\omega_-^2>0$ and
$\omega_-\perp [C,F_1,F_2]$. Since $[C,F_1,F_2]\subset L_\varphi$,
it follows that $\omega_-\in
([C,F_1,F_2])^\perp_{L_\varphi}\otimes\br$.

For $\psi$ we have
$$
L^\psi\otimes \br=L^\varphi\otimes \br\oplus \br
(cl(F_1)-cl(F_2))\oplus \br cl(F)
$$
and
$$
L_\psi\otimes \br= [C,F_1,F_2]^\perp_{L_\varphi}\otimes \br\oplus
\br h\,. $$ Their orthogonal complements to $S=[C,F]$ consist of
$L^\varphi\otimes \br\oplus \br(cl(F_1)-cl(F_2))$ and
$[C,F_1,F_2]^\perp_{L_\varphi}\otimes \br$ respectively. For
the ellipsoid deformation with the action $(L,\tau,\psi)$, we must take
for periods of real K3
$$
\widetilde{\omega}_+\in \left(L^\varphi\otimes \br\oplus \br
(cl(F_1)-cl(F_2))\right)-L^\varphi\otimes \br
$$
instead of $\omega_+\in L^\varphi\otimes \br$
(here $L^\varphi\otimes \br\subset \left(L^\varphi\otimes
\br\oplus \br (cl(F_1)-cl(F_2))\right)$ has codimension one), and
$$
\widetilde{\omega}_-\in [C,F_1,F_2]^\perp_{L_\varphi}\otimes \br
$$
which is the same space as for $\omega_-\in
[C,F_1,F_2]^\perp_{L_\varphi}\otimes \br$. They can be taken
closed to $\omega_+$ and $\omega_-$ and correspond to periods of a
real K3 surface with non-symplectic involution of ellipsoid type
by Global Torelli Theorem for K3 surfaces \cite{PS71}.

Thus, for ellipsoid deformation, we must just replace
$(L,\tau,\varphi)$ by $(L,\tau,\psi)$ where $\psi$ is given by
\eqref{psiF2}. For $(L,\tau,\wvarphi)$ we must replace $\psi$ by
the related involution $\wpsi=\tau\psi$ which is given by the same
formula \eqref{psiF2}.

In \cite{NikulinSaito05}, for an ellipsoid involution
$(L,\tau,\psi)$ the same invariants
$(r=r(\psi),a=a(\psi),\delta=\delta_{\psi S})$ were used (the
invariant $v=v(\psi)$ follows from these invariants). They define
the connected component of moduli of a real K3 involution of
ellipsoid type and corresponding positive real non-singular
bidegree $(4,4)$ curve on ellipsoid. It is easy to calculate these
invariants. The main formula is
\begin{equation}
r(\wpsi)=r(\wvarphi),\ \ a(\wpsi)=a(\wvarphi),\ \
\delta_{S\wpsi}=\delta_{S\wvarphi } \label{wphiwpsiF2}
\end{equation}
(it is geometrically obvious). The remaining relations follow from
\begin{equation}
r(\psi)+r(\wpsi)=22,\ \ a(\psi)=a(\wpsi),\ \delta_{\psi
S}=\delta_{\wpsi S} \label{relatedellF2}
\end{equation}
(see formulae (5.2) and (5.3) in \cite{NikulinSaito05}) and
\eqref{relatedF2}. We obtain
\begin{equation}
r(\psi)=r(\varphi)+2,\ \ a(\psi)=a(\varphi)+2\delta_F,
\delta_{\psi S}=\delta_{\varphi S},\ v(\psi)=v(\varphi)+cl(F)\mod
2. \label{phipsiF2}
\end{equation}
Here \eqref{phipsiF2} follows from Theorem 11 from
\cite{NikulinSaito05} applied to $\tau=-s_{[F_1,F_2]}$ and
$\varphi$, which gives the formal proof of formulae
\eqref{phipsiF2} and \eqref{wphiwpsiF2}.

Thus, the ellipsoid deformation can be represented symbolically as
follows:
\begin{equation}
\begin{array}{ll}
(\bff_2:r(\varphi),a(\varphi),\delta_F(\varphi)=0,\delta_{\varphi
S},v(\varphi)) &\Longrightarrow
(\bee:r=r(\varphi)+2,a=a(\varphi),\delta=\delta_{\varphi S});\\
(\bff_2:r(\varphi),a(\varphi),\delta_F(\varphi)=1,\delta_{\varphi
S},v(\varphi)) &\Longrightarrow
(\bee:r=r(\varphi)+2,a=a(\varphi)+2,\delta=\delta_{\varphi S});\\
(\bff_2:r(\wvarphi),a(\wvarphi),H(\wvarphi)=[F],\delta_{\wvarphi
S},v(\wvarphi)) &\Longrightarrow
(\bee:r=r(\wvarphi),a=a(\wvarphi),\delta=\delta_{\wvarphi S});\\
(\bff_2:r(\wvarphi),a(\wvarphi),H(\wvarphi)=[S],\delta_{\wvarphi
S},v(\wvarphi)) &\Longrightarrow
(\bee:r=r(\wvarphi),a=a(\wvarphi),\delta=\delta_{\wvarphi S}).
\label{elldeformF2}
\end{array}
\end{equation}
There are four types of deformations.

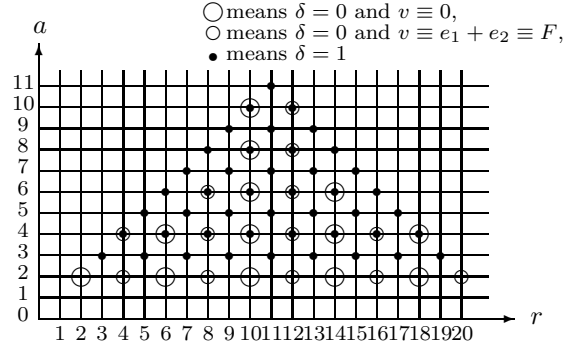
\begin{figure}
\begin{picture}(200,140)
\put(66,116){\circle{7}} \put(71,114){{\tiny means $\delta=0$ and
$v\equiv 0$,}} \put(66,108){\circle{5}} \put(71,106){{\tiny means
$\delta=0$ and $v\equiv e_1+e_2\equiv F$,}} \put(66,
99){\circle*{3}} \put(71, 98){{\tiny means $\delta=1$}}
\multiput(8,0)(8,0){20}{\line(0,1){94}}
\multiput(0,8)(0,8){11}{\line(1,0){170}}
\put(0,0){\vector(0,1){104}} \put(0,0){\vector(1,0){180}} \put(
6,-8){{\tiny $1$}} \put( 14,-8){{\tiny $2$}} \put( 22,-8){{\tiny
$3$}} \put( 30,-8){{\tiny $4$}} \put( 38,-8){{\tiny $5$}} \put(
46,-8){{\tiny $6$}} \put( 54,-8){{\tiny $7$}} \put( 62,-8){{\tiny
$8$}} \put( 70,-8){{\tiny $9$}} \put( 76,-8){{\tiny $10$}} \put(
84,-8){{\tiny $11$}} \put( 92,-8){{\tiny $12$}}
\put(100,-8){{\tiny $13$}} \put(108,-8){{\tiny $14$}}
\put(116,-8){{\tiny $15$}} \put(124,-8){{\tiny $16$}}
\put(132,-8){{\tiny $17$}} \put(140,-8){{\tiny $18$}}
\put(148,-8){{\tiny $19$}} \put(156,-8){{\tiny $20$}}

\put(-8, -1){{\tiny $0$}} \put(-8,  7){{\tiny $1$}} \put(-8,
15){{\tiny $2$}} \put(-8, 23){{\tiny $3$}} \put(-8, 31){{\tiny
$4$}} \put(-8, 39){{\tiny $5$}} \put(-8, 47){{\tiny $6$}} \put(-8,
55){{\tiny $7$}} \put(-8, 63){{\tiny $8$}} \put(-8, 71){{\tiny
$9$}} \put(-10, 79){{\tiny $10$}} \put(-10, 87){{\tiny $11$}}

\put( -2,108){{\footnotesize $a$}} 
\put(186, -2){{\footnotesize $r$}} 

\put( 16,16){\circle{7}} \put( 48,16){\circle{7}} \put(
80,16){\circle{7}} \put(112,16){\circle{7}}
\put(144,16){\circle{7}} \put( 48,32){\circle{7}} \put(
80,32){\circle{7}} \put(112,32){\circle{7}}
\put(144,32){\circle{7}} \put( 80,48){\circle{7}}
\put(112,48){\circle{7}} \put( 80,64){\circle{7}} \put(
80,80){\circle{7}}

\put( 32,16){\circle{5}} \put( 64,16){\circle{5}} \put(
96,16){\circle{5}} \put(128,16){\circle{5}}
\put(160,16){\circle{5}} \put( 32,32){\circle{5}} \put(
64,32){\circle{5}} \put( 96,32){\circle{5}}
\put(128,32){\circle{5}} \put( 64,48){\circle{5}} \put(
96,48){\circle{5}} \put( 96,64){\circle{5}} \put(
96,80){\circle{5}}

\put( 24,24){\circle*{3}} \put( 40,24){\circle*{3}} \put(
56,24){\circle*{3}} \put( 72,24){\circle*{3}} \put(
88,24){\circle*{3}} \put(104,24){\circle*{3}}
\put(120,24){\circle*{3}} \put(136,24){\circle*{3}}
\put(152,24){\circle*{3}} \put( 32,32){\circle*{3}} \put(
48,32){\circle*{3}} \put( 64,32){\circle*{3}} \put(
80,32){\circle*{3}} \put( 96,32){\circle*{3}}
\put(112,32){\circle*{3}} \put(128,32){\circle*{3}}
\put(144,32){\circle*{3}} \put( 40,40){\circle*{3}} \put(
56,40){\circle*{3}} \put( 72,40){\circle*{3}} \put(
88,40){\circle*{3}} \put(104,40){\circle*{3}}
\put(120,40){\circle*{3}} \put(136,40){\circle*{3}} \put(
48,48){\circle*{3}} \put( 64,48){\circle*{3}} \put(
80,48){\circle*{3}} \put( 96,48){\circle*{3}}
\put(112,48){\circle*{3}} \put(128,48){\circle*{3}} \put(
56,56){\circle*{3}} \put( 72,56){\circle*{3}} \put(
88,56){\circle*{3}} \put(104,56){\circle*{3}}
\put(120,56){\circle*{3}} \put( 64,64){\circle*{3}} \put(
80,64){\circle*{3}} \put( 96,64){\circle*{3}}
\put(112,64){\circle*{3}} \put( 72,72){\circle*{3}} \put(
88,72){\circle*{3}} \put(104,72){\circle*{3}} \put(
80,80){\circle*{3}} \put( 96,80){\circle*{3}} \put(
88,88){\circle*{3}}
\end{picture}

\caption{$\bee$: All possible $(r,a,\delta)$ for ellipsoid.}
\label{ellipsoid-graph}
\end{figure}

All invariants $(r,a,\delta)$ for ellipsoid are presented in
Figure \ref{ellipsoid-graph} (it is Figure 16 from
\cite{NikulinSaito05}).

\medskip

From \eqref{elldeformF2} and Figures \ref{diag0F2}, \ref{diag1F2},
\ref{ellipsoid-graph} and \eqref{relatedF2} we see all
deformations of curves from $\bff_2$ to ellipsoid. It is hard to
formulate similar to Theorem \ref{thedefF2hyperbF2} result because
there are too many cases. We restrict by the following.

\begin{theorem} Any positive real non-singular
bidegree $(4,4)$ curve on ellipsoid can be obtained by ellipsoid
deformation of a positive real non-singular curve in
$|-2K_{\bff_2}|$. There are four isomorphism classes of such
deformations at most. All of them are presented in
\eqref{elldeformF2}, \eqref{relatedF2} and Figures
\ref{ellipsoid-graph}, \ref{diag0F2}, \ref{diag1F2}.
\label{thedefF2ellF2}
\end{theorem}

\subsection{Application: Classification of connected components of moduli of
real hyper-elliptically polarized K3 surfaces.}
\label{subsecconnhyp-ellpolK3F2}

For a positive real non-singular curve $A\in |-2K_Y|$ where
$Y=\bff_2$,  or for the corresponding real K3 surface with
degenerate non-symplectic involution $(X,\tau,\varphi)$ of type
$(2,2,0)$, let us consider $P=n_1C+F$ for $n_1\ge 3$. Here $P$ is
primitive and has the degree $n=P^2=4n_1-4$ where $n\ge 8$ and
$n\equiv 0\mod 4$. Thus, $n_1=(n/4)+1$ and $P=((n/4)+1)C+F$. We
have $C\cdot P=2$ for the elliptic pencil $|C|$ on $X$. Then the
complete linear system $|P|$ is hyper-elliptic and gives a double
covering $|P|:X\to Y=\bff_2\subset \bp^N$ where $N=P^2/2+1=n/2+1$.
Here the linear system $|P|$ is pre-image $\pi^\ast|n_1c+s|$ of
the standard complete linear system $|n_1c+s|$ from $Y$. On the
other hand, $P$ is ample since $|P|$ gives a double covering
ramified in a non-singular curve. Thus, $(X,P)$ is a real
hyper-elliptically polarized K3 surface. See \cite{Saint-Donat74},
\cite{NikulinSaito05} and \cite{Nikulin05} for more details.

Thus, our classification of the connected components of moduli
also enumerates connected components of moduli of real
hyper-elliptically polarized K3 surfaces $(X,P)$ of that type. We
shall call them as {\it hyper-elliptically polarized K3 surfaces
of $\bff_2$ type.}

All other types had been considered and classified in
\cite{NikulinSaito05}. See \cite{NikulinSaito05} and
\cite{Nikulin05} for details. For other types, one has to take
$Y=\bp^2$ or $\bff_m$ with $m=0,\,1,\,4$ instead of $\bff_2$ and
also pick up some standard linear system from $Y$. All these
possibilities were described in Sect. 8 of \cite{NikulinSaito05}
(see also \cite{Nikulin05}). Thus, here we finished classification
of connected components of moduli of real hyper-elliptically
polarized K3 surfaces describing the last $\bff_2$ type which had
not been considered. We remark that for a general
hyper-elliptically polarized K3 surface of $\bff_2$ type the
Picard number is $3$, and then the dimension of moduli is $17$.
For all other types of hyper-elliptically polarized K3 surfaces
the Picard number is $2$ and the dimension of moduli is $18$.

Thus, we obtain the following result.

\begin{theorem}
\label{theK3typeF2} Real hyper-elliptically polarized K3 surfaces
$(X,P)$ of type $\bff_2$ (i. e. when the image $|P|:X\to Y$ is
$Y=\bff_2$) exist only for primitive $P$ and $n=P^2\ge 8$ where
$n\equiv 0\mod 4$. For each such $n$ there are exactly 216
connected components of their moduli (of the dimension 17) which
are enumerated by natural invariants in Theorem
\ref{theconnectedF2}.
\end{theorem}

Similar classification for all other types had been obtained in
\cite{NikulinSaito05} (see also \cite{Nikulin05}). Thus, Theorem
\ref{theK3typeF2} completes this classification.

\subsection{Application: Real polarized K3 surfaces as deformations of
real hyper-elliptically polarized K3 surfaces of type
$\bff_2$.}\label{subsecdefF2K3F2} Let $(X,P)$ be a real
hyper-elliptically polarized K3 surface of type $\bff_2$. Its
small deformation gives then a polarized K3 surface
$(\widetilde{X},\widetilde{P})$ where $\widetilde{P}$ is very
ample and defines a connected component of moduli of real
polarized K3 surfaces which are also deformations of $(X,P)$.
Thus, we ask: What real polarized K3 surfaces are deformations of
real hyper-elliptically polarized K3 surfaces of type $\bff_2$?
See \cite{NikulinSaito05} and \cite{Nikulin05} for further details
about this problem where the same problem had been solved for all
other types of real hyper-elliptically polarized K3 surfaces. Thus
here we finalize these results.

By \cite{Nikulin79}, the connected component of moduli of a real
polarized K3 surface is determined by isomorphism class of the
integral polarized involution $(L,\varphi,P)$ corresponding to
$(X,\tau,\varphi)$ with $P=(n/4+1)cl(C)+cl(F)\in S=L^\tau$ (or the
same for $(X,\tau,\wvarphi)$). Here $n\equiv 0\mod 4$ and $n\ge 4$
(we can consider $n=4$ as well, see below). The main invariants (the genus
invariants) of the integral polarized involution $(L,\varphi,P)$
are the same invariants as we used before:
\begin{equation}
(n;r=r(\varphi),a=a(\varphi),\delta_P;\delta_\varphi,\delta_{\varphi
P}) \label{geninvK3F2}
\end{equation}
Almost in all cases these invariants completely determine the
connected component of moduli of real polarized K3 surfaces (see
\cite{Nikulin79} and \cite{Nikulin05}). Thus, in our case there
are four types of such deformations which can be symbolically
expressed as follows where $n\ge 4$ and $n\equiv 0\mod 4$:
\begin{equation}
\begin{array}{ll}
D_n:(\bff_2;r(\varphi),a(\varphi),\delta_
F(\varphi)=1,\delta_\varphi,\delta_{\varphi F}=1) &\Longrightarrow
(n;r=r(\varphi),a=a(\varphi),\delta_P=1,\delta_\varphi,\delta_{\varphi P}=1);\\
D_n:(\bff_2;r(\varphi),a(\varphi),\delta_
F(\varphi)=0,\delta_\varphi,\delta_{\varphi F}) &\Longrightarrow
(n;r=r(\varphi),a=a(\varphi),\delta_P,\delta_\varphi,\delta_{\varphi
P})\
\text{where}\\
\end{array}
\label{defF2K3F2}
\end{equation}
$$
\delta_P= \left\{\begin{array}{cl}
0 &\mbox{if}\ n\equiv 4\mod 8\\
1 &\mbox{otherwise}
\end{array}\right . ,\
\delta_{\varphi P}=\left\{
\begin{array}{cl}
0 &\mbox{if}\ (n\mod 8,\delta_{\varphi F})=
(4\mod 8,0)\\
1 &\mbox{otherwise}
\end{array}\right .\ ;
$$
$$
\begin{array}{ll}
D_n:(\bff_2;r(\wvarphi),a(\wvarphi),H(\wvarphi)=[S],
\delta_\wvarphi) &\Longrightarrow
(n;r=r(\wvarphi),a=a(\wvarphi),\delta_P=0,\delta_\wvarphi,\delta_{\wvarphi P}=1);\\
D_n:(\bff_2;r(\wvarphi),a(\wvarphi),H(\wvarphi)=[F],
\delta_\wvarphi,\delta_{\wvarphi F}) &\Longrightarrow
(n;r=r(\wvarphi),a=a(\wvarphi),\delta_P,\delta_\wvarphi,\delta_{\wvarphi
P})\
\text{where}\\
\end{array}
$$
$$
\delta_P= \left\{\begin{array}{cl}
0 &\mbox{if}\ n\equiv 4\mod 8\\
1 &\mbox{otherwise}
\end{array}\right . ,\
\delta_{\wvarphi P}=\left\{
\begin{array}{cl}
0 &\mbox{if}\ (n\mod 8,\delta_{\wvarphi F})=
(4\mod 8,0)\\
1 &\mbox{otherwise}
\end{array}\right . .
$$
By \cite{Nikulin79} and \cite{Nikulin05}, the invariants
\eqref{geninvK3F2} determine the connected component of moduli of
real polarized K3 surfaces if $r(\varphi)\le 18$. Thus
\eqref{defF2K3F2} completely determine deformations for
$r(\varphi)\le 18$ and $r(\wvarphi)\le 18$.

By Figures \ref{diag0F2} and \ref{diag1F2}, there are no cases
with $r(\varphi)>18$. By \eqref{relatedF2}, $r(\wvarphi)>18$ only
for one connected component
$(r(\wvarphi)=19,a(\wvarphi)=3,H(\wvarphi)=[S],
\delta_{\wvarphi}=1,\delta_{\wvarphi P}=1)$. Then
$X_\wvarphi(\br)=(T_0)^9$. Thus, the deformation $D_n$ can give
only one connected component of moduli of real polarized K3
surfaces (for each $n\ge 4$, $n\equiv 0\mod 4$). By Theorem 16 in
\cite{Nikulin05}, the connected component of moduli of real
polarized K3 surfaces is {\it the standard connected component.}

The list of all possible invariants \eqref{geninvK3F2} of real
polarized K3 surfaces is known, see \cite{Nikulin79} (and also
\cite{Nikulin05}). Thus, using \eqref{defF2K3F2} and Figures
\ref{diag0F2}, \ref{diag1F2}, and \eqref{relatedF2}, we can find
all real polarized K3 surfaces which are deformations of real
hyper-elliptically polarized K3 surfaces of type $\bff_2$. They
must appear to the right of the deformations \eqref{defF2K3F2}. We
can also find out how many these deformations do exist. Three, at
most. We leave further details to the interested reader.

We obtain the following.

\begin{theorem} All real polarized K3 surfaces which are
deformations of real hyper-elliptically polarized K3 surfaces of
type $\bff_2$ and all isomorphism classes of these deformations are
completely determined by the symbolical deformations
\eqref{defF2K3F2}, and Figures \ref{diag0F2}, \ref{diag1F2} and
\eqref{relatedF2} where $n\equiv 0\mod 4$ and $n\ge 8$. For a
fixed real polarized K3 surface there are three isomorphism
classes of these deformations at most.
\label{thedefF2K3F2}
\end{theorem}

For real hyper-elliptically polarized K3 surfaces of all other
types, similar results were obtained in \cite{NikulinSaito05} and
\cite{Nikulin05}. Thus, Theorem \ref{thedefF2K3F2} finalizes these
results by considering the last possible type.

\medskip

Formally, Theorem \ref{thedefF2K3F2} does not include $n=4$
because $P=C+F$ is not ample for $n=4$. We have $P\cdot F=0$, and
the linear system $|P|$ contracts $F_1+F_2$.

But, considering hyperboloid or ellipsoid deformations
$\bff_2\Longrightarrow \bhh$ (see \eqref{defF2H2}),
$\bff_2\Longrightarrow \bee$ (see \eqref{elldeformF2}), and then
taking $P=2cl(C)+cl(F)=e_1+e_2$, we obtain real hyper-elliptically
polarized K3 surfaces of type $\bff_0$ (over $\br$ they are of
hyperboloid or ellipsoid type). Considering their deformations we
obtain real polarized K3 surfaces with polarization of degree four
(a non-singular quartic in $\bp^3$, in general). These
deformations were studied in \cite{NikulinSaito05}. Taking
composition of these two deformations, we can consider $n=4$ for
$\bff_2$ as well.

The hyperboloid deformation $D_4$ gives the same as $D_4$ in
\eqref{defF2K3F2}, and we obtain:
\begin{equation}
\begin{array}{l}
D_4:(\bff_2;r(\varphi),a(\varphi),\delta_
F(\varphi)=1,\delta_\varphi,\delta_{\varphi F}=1) \stackrel{\bhh}
{\Longrightarrow}
(4;r=r(\varphi),a=a(\varphi),\delta_P=1,\delta_\varphi,\delta_{\varphi P}=1);\\
D_4:(\bff_2;r(\varphi),a(\varphi),\delta_
F(\varphi)=0,\delta_\varphi,\delta_{\varphi F})
\stackrel{\bhh}{\Longrightarrow}
(4;r=r(\varphi),a=a(\varphi),\delta_P=0,\delta_\varphi,\delta_{\varphi
P}=\delta_{\varphi F});\\
D_4:(\bff_2;r(\wvarphi),a(\wvarphi),H(\wvarphi)=[S],
\delta_\wvarphi) \stackrel{\bhh}{\Longrightarrow}
(4;r=r(\wvarphi),a=a(\wvarphi),\delta_P=0,\delta_\wvarphi,\delta_{\wvarphi P}=1);\\
D_4:(\bff_2;r(\wvarphi),a(\wvarphi),H(\wvarphi)=[F],
\delta_\wvarphi,\delta_{\wvarphi F})
\stackrel{\bhh}{\Longrightarrow}
(4;r=r(\wvarphi),a=a(\wvarphi),\delta_P=0,\delta_\wvarphi,\delta_{\wvarphi
P}=\delta_{\wvarphi F}).\\
\end{array}
\label{def4HyF2K3F2}
\end{equation}
Obviously, all deformations $D_n$ in \eqref{defF2K3F2} are
compositions of deformations \eqref{defF2H2} of $\bff_2$ to
hyperboloid, and the deformation $D_n$ of type $(\bhh)_1$ (when
$P=(n/4+1)C+F=(n/4)e_1+e_2$ for hyperboloid) considered in
\cite{NikulinSaito05}. (It had been used in \cite{Nikulin05}.)

For ellipsoid deformation we obtain the following where we denote
by $(L,\psi,P)$ the integral polarized involution of a K3 surface
which is obtained under deformation:
\begin{equation}
\begin{array}{ll}
&D_4:(\bff_2:r(\varphi),a(\varphi),\delta_F(\varphi)=1,\delta_{\varphi
S},v(\varphi)) \stackrel{\bee}\Longrightarrow\\
&\ \ \ \ \ \
(4;r(\psi)=r(\varphi)+2,a(\psi)=a(\varphi)+2,\delta_P=0,\delta_\psi,\delta_{\psi P});\\
&D_4:(\bff_2:r(\varphi),a(\varphi),\delta_F(\varphi)=0,\delta_{\varphi
S},v(\varphi)) \stackrel{\bee}{\Longrightarrow}\\
&\ \ \ \ \ \
(4;r(\psi)=r(\varphi)+2,a(\psi)=a(\varphi),\delta_P=0,\delta_\psi,\delta_{\psi P});\\
&D_4:(\bff_2:r(\wvarphi),a(\wvarphi),H(\wvarphi)=[S],\delta_{\wvarphi
S},v(\wvarphi)) \stackrel{\bee}{\Longrightarrow}\\
&\ \ \ \ \ \
(4;r(\psi)=r(\wvarphi),a(\psi)=a(\wvarphi),\delta_P=0,\delta_\psi,\delta_{\psi
P});\\
&D_4:(\bff_2:r(\wvarphi),a(\wvarphi),H(\wvarphi)=[F],\delta_{\wvarphi
S},v(\wvarphi)) \stackrel{\bee}{\Longrightarrow}\\
&\ \ \ \ \ \
(4;r(\psi)=r(\wvarphi),a(\psi)=a(\wvarphi),\delta_P=0,\delta_\psi,\delta_{\psi P})\\
\end{array}
\label{def4ElF2K3F2}
\end{equation}
where
$$
\delta_\psi= \left\{\begin{array}{cl}
0 &\mbox{if}\ (\delta_{\varphi S}=\delta_{\wvarphi S},\ r(\psi)\mod 4)=(0,2\mod4)\\
1 &\mbox{otherwise}
\end{array}\right . ,
$$
$$
\delta_{\psi P}=\left\{
\begin{array}{cl}
0 &\mbox{if}\ (\delta_{\varphi S}=\delta_{\wvarphi S},\ r(\psi)\mod 4)=(0,0\mod 4)\\
1 &\mbox{otherwise}
\end{array}\right . .
$$
Thus, similarly to Theorem \ref{thedefF2K3F2} we obtain

\begin{theorem} All non-singular real quartics in $\bp^3$ which are
deformations of double coverings of real quadratic cone in $\bp^3$
ramified in a non-singular degree $8$ curve, and all isomorphism
classes of these deformations are completely determined by the
symbolical deformations \eqref{def4HyF2K3F2}, \eqref{def4ElF2K3F2}
and Figures \ref{diag0F2}, \ref{diag1F2} and \eqref{relatedF2}.
For a fixed non-singular real quartic in $\bp^3$ there are seven
isomorphism classes of these deformations at most.
\label{thedefF2quarticF2}
\end{theorem}

Similar results for non-degenerate real quadrics in $\bp^3$ (e. g.
of hyperboloid type, or ellipsoid) had been obtained in
\cite{NikulinSaito05}. Thus, Theorem \ref{thedefF2quarticF2}
completes these results for a non-degenerate quadratic cone in
$\bp^3$.

\bigskip

All our results here (and results of  \cite{NikulinSaito05} and
\cite{Nikulin05} as well) are based on Global Torelli Theorem
\cite{PS71} and surjectivity of Torelli map \cite{Kulikov77} for
K3 surfaces and some geometric and arithmetic considerations.
Thus, our results are mainly existence and uniqueness results. It
would be very interesting to construct exactly all representatives
of connected components of moduli which we considered and
classified. Our results about deformations reduce everything to
construction of representatives of connected components of moduli
of real non-singular curves $A\in |-2K_{\bff_m}|$, $0\le m \le 4$,
or degree 6 curves in $\bp^2$. In many cases they are known by
direct construction. E. g. see \cite{Gudkov69}, \cite{Gudkov79},
\cite{GudkovShustin80}, \cite{Rokhlin78}, \cite{Viro80},
\cite{Wilson}, \cite{Zvonilov82}, \cite{Zvonilov83},
\cite{Zvonilov92}. It would be very interesting to have such
constructions for all connected components of moduli. It will make
our results effective and useful for some questions which don't
depend on deformation. This construction is surely possible
because all our connected components of moduli are exactly
labelled by some natural invariants, and number of the connected
components is bounded.



\end{document}